\input ./style/arxiv-general.cfg
\documentclass[MSNbibl,number,citesort,seceqn,secthm,dvips]{arxbj}
\makeatletter
   \@ifpackageloaded{graphicx}{}{\usepackage{graphicx}}
\makeatother

\volume{22}
\issue{4}
\pubyear{2016}
\firstpage{2401}
\lastpage{2441}
\doi{10.3150/15-BEJ733}
\docsubty{FLA}

\makeatletter
\renewcommand{\bot}{\perp}
\newcommand{\rrvert}{\vert}
\newcommand{\llvert}{\vert}
\newtheorem{theorem}{Theorem}[section]
\newcommand{\eqref}[1]{(\ref{#1})}

\newremark{property}{Property}[section]
\newtheorem{proposition}{Proposition}[section]
\newproclaim{definition}{Definition}[section]
\newremark{remark}{Remark}[section]
\newproclaim{assumption}[theorem]{Assumption}
\def\fddlim{\mathrm{f.d.d.\mbox{-}lim}}
\def\limd{\stackrel{\mathrm{d}}{\longrightarrow}}
\def\limfdd{\stackrel{\mathrm{f.d.d.}}{\longrightarrow}}
\def\fddlimn{\mathop{\mathrm{f.d.d.\mbox{-}lim}}_{n \rightarrow \infty}}
\def\fddlimm{\mathop{\mathrm{f.d.d.\mbox{-}lim}}_{m \rightarrow \infty}}
\def\neqfdd{\stackrel{\mathrm{f.d.d.}}{\neq}}
\def\eqfdd{\stackrel{\mathrm{f.d.d.}}{=}}
\def\eqd{\stackrel{\mathrm{d}}{=}}

\newtheorem{lem}[thm]{Lemma}

\newproclaim{defn}[thm]{Definition}
\newremark{rem}{Remark}[section]
\newremark{ex}[thm]{Example}
\makeatother

\begin{document}
\begin{frontmatter}

\title{Aggregation of autoregressive random fields and anisotropic
long-range dependence}
\runtitle{Aggregation of autoregressive random fields}

\begin{aug}
\author[A]{\inits{D.}\fnms{Donata}~\snm{Puplinskait\normalfont{\.{E}}}\thanksref{A}\ead[label=e1]{donata.puplinskaite@mif.vu.lt}}
\pdfauthor{Donata Puplinskait\.e, Donatas Surgailis}
\and
\author[B]{\inits{D.}\fnms{Donatas}~\snm{Surgailis}\corref{}\thanksref{B}\ead[label=e2]{donatas.surgailis@mii.vu.lt}}
\runauthor{D. Puplinskait\.e and  D. Surgailis} 
\address[A]{Faculty of Mathematics and Informatics, Vilnius University,
Naugarduko 24, LT-03225 Vilnius, Lithuania. \printead{e1}}
\address[B]{Institute of Mathematics and Informatics, Vilnius
University, Akademijos 4, LT-08663 Vilnius, Lithuania. \printead{e2}}
\end{aug}

%
\received{\smonth{12} \syear{2014}}
%
\revised{\smonth{3} \syear{2015}}

\begin{abstract}
We introduce the notions of scaling transition and
distributional
long-range dependence for stationary random fields $Y$ on $\mathbb
{Z}^2$ whose
normalized partial sums on
rectangles with sides growing at rates $O(n) $ and $O(n^{\gamma})$
tend to an operator scaling random field $V_\gamma$ on $\mathbb
{R}^2$, for any
$\gamma>0$. The scaling transition
is characterized by the fact that there exists a unique $\gamma_0 >0$
such that the scaling limits
$V_\gamma$ are different and do not depend on $\gamma$ for $\gamma>
\gamma_0$ and $\gamma< \gamma_0$.
The existence of scaling transition together with
anisotropic and isotropic distributional long-range dependence properties
is demonstrated for a class of $\alpha$-stable $(1 < \alpha\le2)$ aggregated
nearest-neighbor autoregressive random fields on $\mathbb{Z}^2$
with a scalar random coefficient $A$
having a regularly varying probability density near the ``unit root'' $A=1$.
\end{abstract}

\end{frontmatter}

\section{Introduction}

Following Bierm\'e \textit{et al.} \cite{bier2007},
a scalar-valued random field (RF) $ V= \{ V(x); x \in\mathbb{R}^\nu\}
$ is
called \textit{operator scaling random field} (OSRF) if there exist a $H
>0$ and
a $\nu\times\nu$ real matrix $E $ whose all eigenvalues have positive
real parts, such that for any $\lambda>0$
%
\begin{equation}
\label{OS} \bigl\{ V\bigl(\lambda^E x\bigr); x \in
\mathbb{R}^\nu\bigr\} \eqfdd \bigl\{ \lambda^H V(x); x \in
\mathbb{R}^\nu\bigr\}.
\end{equation}
(See the end of this section for all unexplained notation.)
In the case when $E = I$ is the unit matrix, \eqref{OS} agrees with the
definition of $H$-self-similar random field (SSRF), the latter
referred to as self-similar process when $\nu=1$.
OSRFs may exhibit strong anisotropy and play an important
role in various physical theories; see \cite{bier2007} and the
references therein. Several classes of OSRFs were constructed and discussed
in \cite{bier2007,clau2009}.

It is well known that the class of self-similar processes is very
large, SSRFs and OSFRs being even more numerous.
According to a popular view, the ``value'' of a concrete self-similar
process depends on its ``domain of attraction''. In the case $\nu=1$,
the domain of attraction of a self-similar stationary increment process
$V= \{V(\tau); \tau\ge0\}$ is
defined in \cite{lamp1962} as the class of all stationary processes $
Y= \{ Y(t); t \in\mathbb{Z}_+ \}$ whose normalized
partial sums tend to $V$ in the distributional sense, namely,
%
\begin{equation}
\label{Xsum0} B^{-1}_n \sum
_{t=1}^{[n\tau]} Y(t) \limfdd V(\tau), \qquad\tau \in
\mathbb{R}_+, n \to\infty.
\end{equation}
The classical Lamperti's theorem \cite{lamp1962} says that in the case
of \eqref{Xsum0}, the normalizing constants $B_n$ necessarily grow as
$n^H$ (modulus a slowly varying factor) and the limit random process
in~\eqref{Xsum0} is $H$-self-similar. The limit process
$V$ in \eqref{Xsum0} characterizes large-scale and dependence
properties of $Y$, leading to the important
concept of \textit{distributional short/long memory} originating in Cox
\cite{cox1984}; see also
(\cite{deh2002}, pages 76--77),
\cite{gir2009,ps2009,ps2010,pps2013,pils2014}.
There exists a large probabilistic
literature devoted to studying the partial sums limits
of various classes of strongly and weakly dependent processes and RFs.
In particular, several works \cite{dob1979,dobmaj1979,leo1999,sur1982,douk2002,lav2007} discussed
the partial sums limits of (stationary) RFs indexed by $t \in\mathbb
{Z}^\nu$:
%
\begin{equation}
\label{Xsum00} B^{-1}_n \sum
_{t \in K_{[nx]}} Y(t) \limfdd V(x), \qquad x =(x_1, \ldots,
x_\nu) \in\mathbb{R}^\nu_+, n \to\infty,
\end{equation}
where $K_{[nx]} := \{ t = (t_1, \ldots, t_\nu) \in\mathbb{Z}^\nu:
1\le t_i
\le
nx_i \}$ is a sequence of rectangles whose all sides increase as $O(n)$.
Related results for Gaussian or linear (shot-noise) and their
subordinated RFs,
with a particular focus on
large-time behavior of statistical solutions of partial differential equations,
were obtained in~\cite{alb1994,anh1999,leool2012,leo1999,leo2011}. See also the recent
paper Anh \textit{et al.}~\cite{AnhLRM2012} and the numerous references therein.
Most of the above mentioned studies deal with ``nearly isotropic''
models of RFs characterized by a single memory parameter $H$ and a
limiting SSRF $\{V(x)\}$ in \eqref{Xsum00}.

Similarly as in the case of random processes indexed by $\mathbb{Z}$,
stationary RFs usually exhibit two types of dependence: weak dependence
and strong dependence. The second type
of dependence is often called long memory or \textit{long-range
dependence} (LRD). Although there is no single satisfactory
definition of LRD, usually it refers to a stationary RF $Y$ having an
unbounded spectral density $f$:
$\sup_{x \in[-\pi, \pi]^\nu} f(x) = \infty$
or a non-summable auto-covariance function: $\sum_{t \in\mathbb
{Z}^\nu}
| \operatorname{cov}(Y(0), Y(t))| = \infty$; see
\cite{dobmaj1979,douk2002,douk2003p,douk2003,lav2007,book2012,ber2013}. The above definitions of LRD do not apply to
RFs with infinite variance and are of limited use since these
properties are very hard to test
in practice. On the other hand, the characterization of LRD based on
partial sums as in the case of distributional long memory
is directly related to the asymptotic
distribution of the sample mean. As noted in \cite{guo2009}, in many
applications
the auto-covariance of RF decays
with different exponents (Hurst indices) in different directions. In
the latter case,
the partial sums of such RF on rectangles $\prod_{i=1}^\nu[1, n_i] $
may grow at different rate with
$n_i \to\infty$, leading to a limiting anisotropic OSRF.

The present paper attempts a systematic study of \textit{anisotropic
distributional long-range dependence}, by exhibiting some natural classes
of RFs whose partial sums tend to OSRFs.
Our study is limited to the case $\nu=2$ and
RFs with anisotropy along the coordinate axes and a diagonal matrix $E$.
Note that for $\nu=2$ and
$E = \operatorname{diag}(1, \gamma),  0< \gamma\ne1$, relation (\ref{OS})
writes as
$\{ V(\lambda x, \lambda^\gamma y) \} \eqfdd\{ \lambda^H V(x,y) \}$, or
%
\begin{equation}
\label{OS1} \bigl\{ \lambda V(x,y); (x,y) \in\mathbb{R}^2 \bigr\}
\eqfdd \bigl\{ V\bigl(\lambda ^{1/H} x, \lambda^{\gamma/H} y\bigr);
(x,y) \in\mathbb{R}^2 \bigr\} \qquad \forall \lambda> 0.
\end{equation}
The OSRFs $V = V_\gamma$
depending on $\gamma>0$ are obtained by taking the partial sums limits
%
\begin{equation}
\label{Xsum01} n^{-H(\gamma)} \sum_{(t,s) \in K_{[nx, n^{\gamma}y]}} Y(t,s)
\limfdd V_\gamma(x,y), \qquad (x,y) \in\mathbb{R}^2_+, n \to
\infty
\end{equation}
on rectangles
$K_{[nx, n^{\gamma}y]} := \{ (t,s) \in\mathbb{Z}^2: 1\le t \le nx,
1\le s
\le
n^{\gamma}y \} $ whose sides grow at possibly different
rate $O(n)$ and $O(n^{\gamma})$. Somewhat unexpectedly, it turned out
that for a large class of RFs
$ Y= \{Y(t,s); (t,s)\in\mathbb{Z}^2 \}$, the limit in \eqref{Xsum01} exists
for any $\gamma>0$.
What is more surprising, many LRD RFs $Y$
exhibit a dramatic change of their scaling behavior
at some point $\gamma_0 >0$, in the sense that $V_\gamma\eqfdd V_\pm$
do not depend on
$\gamma$ for $ \gamma> \gamma_0 $ or $\gamma< \gamma_0$ and $V_+
\neqfdd V_-$.
This
phenomenon which we call
\textit{scaling transition} seems to be of general nature, suggesting an
exciting new area in spatial research \cite{ps2014}. It
occurs for $\alpha$-stable ($1< \alpha\le2$) aggregated
autoregressive RFs
studied in this paper, for a natural class of LRD Gaussian RFs
discussed in \cite{ps2014} and Remark~\ref{remGaus} below, but also in
a very different
context of network traffic models; see
Remark~\ref{remON}. In most of the above
mentioned works, the limit $V_{\gamma_0}$ is different from $V_+ $ and
$V_-$, and the differences
between $V_{\gamma_0}, V_+, V_-$ can be characterized by dependence properties
of increments $V_\gamma(K) := V_\gamma(x,y) - V_\gamma(u,y) -
V_\gamma
(x,v) + V_\gamma(u,v) $ on rectangles $K = (u,x] \times(v,y] \subset
\mathbb{R}^2_+$,
which may change from independent increments in the vertical
direction for $\gamma> \gamma_0$
to independent increments in the horizontal direction (or completely
dependent increments in the vertical direction)
for $\gamma< \gamma_0$, or vice versa.
Further on, depending on whether $\gamma_0 = 1 $ or $\gamma_0 \ne1$,
the corresponding RF $Y$ is said to have \textit{isotropic distributional
LRD} or
\textit{anisotropic distributional LRD} properties.

\begin{figure}[b]
\begin{tabular}{@{}cc@{}}

\includegraphics{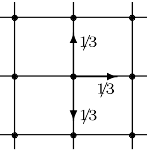}
 & \includegraphics{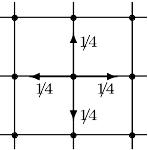}\\
\footnotesize{(a) 3N} & \footnotesize{(b) 4N}
\end{tabular}
\caption{One-step transition probabilities of the random walk underlying models \eqref{3N}
and \eqref{4N}.}\label{fig1}
\vspace*{-6pt}
\end{figure}

The main purpose of this work is establishing scaling transition and
Type I isotropic and anisotropic distributional LRD properties for
a natural class of aggregated nearest-neighbor random-coefficient
autoregressive RFs with finite and infinite variance. We recall that
the idea of contemporaneous aggregation originates to Granger \cite
{gran1980}, who
observed that aggregation of random-coefficient AR(1) equations with
random beta-distributed coefficient can lead to a Gaussian
process with long memory and slowly decaying covariance function. Since
then, aggregation became one
of the most important methods for modeling and studying long memory
processes; see Beran \cite{ber2013}.
For linear and heteroscedastic autoregressive time series models with
one-dimensional time, it was developed in
\cite{gon1988,opp2004,kaz2004,zaff2004,zaff2007,cel2007,gir2010,ps2009,ps2010,pps2013,pils2014} and for some RF
models in \cite{lav2007,lav2011,lls2014,leot2013}.
Aggregation is also important for understanding and modeling of spatial
LRD processes
by relating them to short-range dependent random-coefficient
autoregressive models in a natural way.
The two models of interest are given by equations:
%
\begin{eqnarray}
\label{3N} X_3(t,s) &=&\frac{A}{3} \bigl(X_3(t-1,s)
+ X_3(t,s+1) + X_3(t,s-1) \bigr) + \varepsilon (t,s),
\\
\label{4N} X_4(t,s)&=&\frac{A}{4} \bigl(X_4(t-1,s)
+ X_4(t+1,s) + X_4(t,s+1) + X_4(t,s-1)
\bigr) + \varepsilon (t,s),\quad
\end{eqnarray}
where $\{\varepsilon (t,s), (t,s) \in\mathbb{Z}^2 \} $ are i.i.d.
r.v.'s whose generic
distribution $\varepsilon $ belongs
to the domain of (normal) attraction of $\alpha$-stable law,
$1< \alpha\le2$, and
$A \in[0,1)$ is a r.v. independent of $\{\varepsilon (t,s)\} $ and
having a regularly varying probability density $\phi$ at $a=1$: there
exist $\phi_1 >0$ and $\beta> -1 $ such that
%
\begin{eqnarray}
\label{mixdensity} \phi(a)&\sim&\phi_1 (1-a)^\beta, \qquad a
\nearrow1.
\end{eqnarray}
In the sequel, we refer to \eqref{3N} and \eqref{4N} as the 3N and 4N models,
N standing for ``Neighbor''. Let $X_{3j}, X_{4j}, j=1, \ldots, m$ denote
$m$ independent copies of $X_3, X_4$ in \eqref{3N}, \eqref{4N},
respectively. As shown in Section~\ref{sec5}, the aggregated
3N and 4N models defined as $m^{-1} \sum_{j=1}^m X_{ij}(t,s) \limfdd
\mathfrak{X}_i (t,s), m \to\infty, i=3,4$
are written as respective mixed $\alpha$-stable moving-averages:
%
\begin{eqnarray}
\label{aggre34} \mathfrak{X}_i(t,s) &=&\sum
_{(u,v) \in\mathbb{Z}^2} \int_{[0,1)} g_i(t-u,s-v,
a) M_{u,v}(\mathrm{d}a), \qquad (t,s) \in\mathbb{Z}^2,
i=3,4,
\end{eqnarray}
where $\{ M_{u,v}(\mathrm{d}a),  (u,v) \in\mathbb{Z}^2 \} $ are
i.i.d. copies of an
$\alpha$-stable random measure $M$ on $[0,1)$ with control measure
$\phi(a) \mathrm{d}a$
and $g_i$ is the corresponding (lattice) Green function:
%
\begin{eqnarray}
\label{green00} g_i(t,s,a)&=&\sum_{k=0}^\infty
a^{k} p_k(t,s), \qquad (t,s) \in\mathbb{Z}^2,
a \in[0,1), i=3,4,
\end{eqnarray}
where $p_k(t,s) = \mathrm{P}(W_k = (t,s) |W_0 = (0,0))$ is the $k$-step
probability of the nearest-neighbor random walk $\{ W_k, k =0,1,\ldots
\} $
on the lattice $\mathbb{Z}^2$ with one-step transition
probabilities $p(t,s)$ shown in Figure~\ref{fig1}(a)--(b).

The main results of Sections~\ref{sec3} and \ref{sec4} are Theorems \ref{3Ntheo} and
\ref
{4Ntheo}. The first theorem identifies the scaling limits $V_\gamma,
\gamma>0 $ in \eqref{Xsum01} and proves
Type I anisotropic LRD
property in the sense of Definition~\ref{Adlm} with $\gamma_0 = 1/2 $
for the aggregated 3N model $\mathfrak{X}_3$ in \eqref{aggre34}.
Similarly, the second theorem obtains Type I isotropic LRD property
($\gamma_0 = 1$) for the aggregated 4N model $\mathfrak{X}_4$ in~\eqref{aggre34}.

The proofs of Theorems \ref{3Ntheo} and \ref{4Ntheo} rely
on the asymptotics of the lattice Green function in~\eqref{green00}
for models
3N and 4N. Particularly, Lemmas \ref{lemma3N} and \ref{lemma4N} obtain
the following point-wise convergences: as $\lambda\to\infty$,
%
\begin{eqnarray}
\label{green03} \sqrt{\lambda} g_3 \biggl([\lambda t], [\sqrt{
\lambda}s], 1- \frac
{z}{\lambda} \biggr) &\to& h_3(t,s,z), \qquad t
>0, s \in\mathbb{R}, z>0,
\\
\label{green04} g_4 \biggl([\lambda t], [\lambda s], 1-
\frac{z}{\lambda^2} \biggr) &\to&h_4(t,s,z), \qquad (t,s) \in
\mathbb{R}^2_0, z>0,
\end{eqnarray}
respectively,
together with dominating bounds of the left-hand sides of \eqref{green03}--\eqref{green04}.
The limit functions $h_3$ and $h_4 $ in \eqref{green03}--\eqref{green04}
(entering stochastic integral representations
of the scaling limits $V_\gamma$ in Theorems \ref{3Ntheo} and \ref{4Ntheo})
are given by
%
\begin{eqnarray}
h_3(t,s,z) &:=& \frac{3}{2\sqrt{\pi t}}\mathrm{e}^{-3zt - {s^2}/({4t})}
\mathbf{1} (t,z >0),
\nonumber
\\[-8pt]
\label{hgreen}
\\[-8pt]
\nonumber
h_4(t,s,z) & :=& \frac{2}{\pi}K_0
\bigl(2\sqrt{z\bigl(t^2 + s^2\bigr)} \bigr)
\mathbf{1}(z>0),
\end{eqnarray}
where $K_0$ is the modified Bessel function of second kind. Note that
$h_3$ in \eqref{hgreen} is the Green function of one-dimensional heat
equation (modulus constant coefficients), while $h_4$ is the Green
function of the Helmholtz equation in $\mathbb{R}^2$.
The proofs of these technical lemmas can be found in the extended
version of this paper available at \surl{http://arxiv.org/abs/1303.2209v3} and
will be published elsewhere.
Lemmas \ref{lemma3N} and \ref{lemma4N}
may also have independent interest for studying the behavior
of the autoregressive fields \eqref{3N} and \eqref{4N} with
deterministic coefficient $A$ in the vicinity
of $A=1$, particularly,
for testing stationarity near the unit root in spatial autoregressive models,
cf. \cite{bhat1997}.

\vspace*{6pt}

\textit{Notation}. \
In what follows, $C, C(K), \ldots$ denote generic
constants, possibly depending on the variables in brackets,
which may be different at different locations. We write $\limd$, $\eqd$,
$\limfdd$, $\eqfdd$, $\neqfdd$ for the weak convergence and equality
and inequality of distributions and finite-dimensional
distributions, respectively. $\fddlim$ stands for
the limit in the sense of weak convergence of finite-dimensional distributions.
For $\lambda>0$ and a $\nu\times\nu$ matrix $E,\lambda^E :=
\mathrm{e}^{E
\log\lambda}$, where
$\mathrm{e}^A = \sum_{k=0}^\infty A^k/k! $ is the matrix exponential.
$\mathbb{Z}^\nu_+
:= \{ (t_1, \ldots, t_\nu) \in\mathbb{Z}^\nu: t_i > 0,
i=1, \ldots, \nu\}, \mathbb{R}^\nu_+ := \{ (x_1, \ldots, x_d) \in
\mathbb{R}^\nu
: x_i
> 0,
i=1, \ldots, \nu\}, \bar\mathbb{R}^\nu_+ := \{ (x_1, \ldots, x_d)
\in\mathbb{R}
^\nu:
x_i \ge0,
i=1, \ldots, \nu\},
\mathbb{Z}_+ := \mathbb{Z}_+^1, \mathbb{R}_+ := \mathbb{R}^1_+,
\bar\mathbb{R}_+ := \bar\mathbb{R}^1_+,  \mathbb{R}
^2_0 :=
\mathbb{R}^2 \setminus\{(0,0)\}$. $E = \operatorname{diag}(\gamma_1, \ldots,
\gamma
_\nu
)$ denotes
the diagonal $\nu\times\nu$ matrix with entries $\gamma_1, \ldots,
\gamma_\nu$ on the diagonal. $\mathbf{1}_A $ stands for the indicator
function of
a set $A$. $\log_+(x) := \log x, x \ge1, := 0 $ otherwise. $[x] =
\lfloor
x \rfloor := k,  x \in[k, k+1),
\lceil y \rceil := k+ 1, y \in(k, k+1],  k \in\mathbb{Z}$.
$K_{[nx, n^{\gamma}y]}
:= \{ (t,s) \in\mathbb{Z}^2: 1\le t \le nx, 1\le s \le n^{\gamma}y
\} $,
$K_{(u,v); (x,y)} := \{ (t,s) \in\mathbb{R}^2_+: u < t \le x, v < s
\le y \} $.

\section{Scaling transition and Type I distributional LRD for RFs on~$\mathbb{Z}^2$}\label{sec2}

In this section, by RF on $\bar\mathbb{R}^2_+$ we mean a RF $V = \{V(x,y);
(x,y) \in\bar\mathbb{R}^2_+\}$ such that
$V(x,y) = 0$ for any $(x,y) \in\bar\mathbb{R}^2_+\setminus\mathbb{R}^2_+$.
A RF $V$ on $\bar\mathbb{R}^2_+$ is said \textit{trivial} if $V(x,y) =
0$ for any
$(x,y) \in\bar\mathbb{R}^2_+$, else $V$ is said \textit{non-trivial}.

\begin{defn} \label{phase}
Let\vspace*{1pt} $ Y = \{Y(t,s); (t,s) \in\mathbb{Z}^2
\}$ be
a stationary RF. Assume that for any $\gamma>0$ there exist a
normalization $A_n(\gamma) \to\infty$
and a non-trivial
RF  $V_\gamma= \{V_\gamma(x,y); (x,y) \in\bar\mathbb{R}^2_+\}$
such that
%
\begin{equation}
\label{Xsum} A^{-1}_n(\gamma) \sum
_{(t,s) \in K_{[nx, n^{\gamma}y]}} Y(t,s) \limfdd V_\gamma(x,y), \qquad (x,y) \in
\mathbb{R}^2_+, n \to \infty.
\end{equation}
We say that  $Y$ \textup{exhibits scaling transition} if there exists
$\gamma_0> 0$ such that
the limits $V_\gamma\eqfdd V_+, \gamma> \gamma_0$
and $V_\gamma\eqfdd V_-,  \gamma<\gamma_0 $ do not depend on
$\gamma$ for
$\gamma>\gamma_0$ and $\gamma< \gamma_0$
and, moreover, $V_+$ and $V_- $ are mutually different RFs, in the
sense that
for any $a >0$\vspace*{-3pt}
%
\begin{equation}
\label{ac} V_+ \neqfdd aV_-.
\end{equation}
In such case, $V_{\gamma_0}$ will be called the \textup{well balanced} and
$V_+, V_-$
the \textup{unbalanced} scaling
limits of~$Y$, respectively.
\end{defn}

Note that the fact that \eqref{ac} hold for any $a >0$
excludes a trivial change of the scaling limit by
a linear change of normalization. It follows rather easily that under
general set-up scaling limits $V_\gamma$ satisfy the
self-similarity and stationarity of rectangular increments properties
stated in Proposition~\ref{Vanis} below.
Let $V = \{V(x,y); (x,y) \in\bar\mathbb{R}^2_+ \}$ be a
RF and $K = K_{(u,v); (x,y)} \subset\mathbb{R}^2_+ $ be a rectangle.
By \textit{increment of $V$ on rectangle $K$} we mean the\vspace*{-2pt} difference\looseness=-1
\[
V(K) := V(x,y) - V(u,y) - V(x,v) + V(u,v).
\]
We say that $V$ \textit{has stationary rectangular increments} if for any
$(u,v) \in\mathbb{R}^2_+$,\vspace*{-2pt}
%
\begin{equation}
\label{Vstatinc} \bigl\{ V(K_ {(u,v); (x,y)}); x\ge u, y \ge v \bigr\} \eqfdd\bigl
\{V(K_ {(0,0); (x-u,y-v)}); x\ge u, y \ge v \bigr\}.
\end{equation}

\begin{figure}[t]

\includegraphics{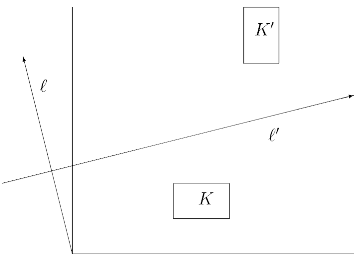}

\caption{Rectangles $K$ and $K'$ separated by line $\ell'$.}\label{fig2}
\end{figure}

As mentioned in the Introduction,
in the case of scaling transition the limits $V_{\gamma_0}, V_+, V_-$
can be characterized by
dependence properties of increments $V(K)$. To define these properties,
we introduce some terminology.
Let $\ell= \{ (x,y) \in\mathbb{R}^2: a x + b y = c \} $ be a line in
$\mathbb{R}^2$.
A line
$\ell' = \{ (x,y) \in\mathbb{R}^2: a' x + b' y = c' \} $ is said
\textit{perpendicular to} $\ell$  (denoted $\ell' \bot\ell$) if $a a'+
b b' = 0$.
We say that two
rectangles $K = K_{(u,v); (x,y)} $ and $K'= K_{(u',v'); (x',y')} $
are \textit{separated by line $\ell' $} if they lie on different sides
of $\ell' $, in which case $K$ and $K'$ are necessarily disjoint: $K
\cap K' = \varnothing$.
See Figure~\ref{fig2}.

%

\begin{defn} \label{iinc}
Let $V= \{V(x,y); (x,y) \in\bar\mathbb{R}^2_+ \}$ be a RF with stationary
rectangular increments,
$V(x,0) = V(0,y) \equiv0,  x,y \ge0$, and $\ell\subset\mathbb
{R}^2$ be a
given line , $(0,0) \in\ell$.
We say that $V$ has:
\begin{longlist}[(iii)]
\item[(i)]  \textup{independent\vspace*{1pt} rectangular
increments in direction} $\ell$   if for any orthogonal line $\ell' \bot\ell$
and any two rectangles $K, K' \subset\mathbb{R}^2_+$ separated by
$\ell'$, increments
$V(K)$ and $ V(K')$
are independent;

\item[(ii)] \textup{invariant rectangular
increments
in direction $\ell$} if $V(K) = V(K')$
for any two rectangles $K, K' \subset\mathbb{R}^2_+$ such that
$K' = (x,y) + K$ for some $(x,y) \in\ell$;

\item[(iii)]  \textup{properly dependent rectangular
increments
in direction $\ell$} if neither \textup{(i)} nor \textup{(ii)} holds;

\item[(iv)] \textup{properly dependent rectangular
increments} if $V$ has properly dependent rectangular increments
in arbitrary direction;

\item[(v)] \textup{independent rectangular
increments} if $V$ has independent rectangular increments
in arbitrary direction.
\end{longlist}
\end{defn}

\begin{ex} \label{exFBS}
\textit{Fractional Brownian sheet} $B_{H_1,H_2}$ with parameters $0 <
H_1, H_2 \le1$ is a
Gaussian process on $\bar\mathbb{R}^2_+$ with zero mean and covariance
%
\begin{eqnarray}
&& \mathrm{E}B_{H_1,H_2}(x,y) B_{H_1,H_2}
\bigl(x',y'\bigr)
\nonumber
\\[-8pt]
\label{covB}
\\[-8pt]
\nonumber
&&\quad= \frac{1}{4}
\bigl(x^{2H_1} + x'^{2H_1} - \bigl|x-x'\bigr|^{2H_1}
\bigr) \bigl(y^{2H_2} + y'^{2H_2} -
\bigl|y-y'\bigr|^{2H_2}\bigr),
\end{eqnarray}
where $(x,y), (x',y') \in\bar\mathbb{R}^2_+$. It follows (see \cite{aya2002},
Corollary~3)
that for any rectangles
$K = K_{(u,v); (x,y)}, K' = K_{(u',v'); (x',y')}$
%
\begin{eqnarray}
&&\mathrm{E}B_{H_1,H_2}(K) B_{H_1,H_2}\bigl(K'\bigr)
\nonumber
\\
&&\quad = \tfrac{1}{4} \bigl(\bigl|x-x'\bigr|^{2H_1} +
\bigl|u-u'\bigr|^{2H_1} - \bigl|x-u'\bigr|^{2H_1} -
\bigl|x'-u\bigr|^{2H_1} \bigr)
\nonumber
\\
\label{FBsheet} && \qquad{}\times \bigl(\bigl|y-y'\bigr|^{2H_2} +
\bigl|v-v'\bigr|^{2H_2} - \bigl|y-v'\bigr|^{2H_2} -
\bigl|y'-v\bigr|^{2H_2} \bigr)
\\
&&\qquad= \mathrm{E}\bigl(B_{H_1}(x)-B_{H_1}(u)\bigr)
\bigl(B_{H_1}\bigl(x'\bigr)-B_{H_1}
\bigl(u'\bigr)\bigr)
\nonumber
\\
&&\qquad\quad{}\times\mathrm {E}\bigl(B_{H_2}(y)-B_{H_2}(v)
\bigr) \bigl(B_{H_2}\bigl(y'\bigr)-B_{H_2}
\bigl(v'\bigr)\bigr),\nonumber
\end{eqnarray}
where $\{B_{H}(x); x\in\bar\mathbb{R}_+\}$ is a fractional Brownian
motion on
$\bar\mathbb{R}_+ = [0,\infty)$ with
$\mathrm{E}B_{H}(x)\* B_{H}(x') = (1/2)(x^{2H} + x'^{2H} - |x-x'|^{2H}),
 H \in
(0,1]$. (Recall that
$B_{1/2}$ is a standard Brownian motion with variance $\mathrm
{E}B^2_{1/2} (x)
= x $ and
$B_{1}(x) = x B_1(1)$ is a random line.)
In particular, $B_{H_1,H_2}$ has stationary rectangular increments; see
\cite{aya2002}, Proposition~2. It follows from \eqref{FBsheet} that
$B_{1/2,H_2}$ has independent rectangular increments in the horizontal
direction since
$\mathrm{E}B_{1/2,H_2}(K) B_{1/2,H_2}(K') = 0$ for any $K, K'$ which are
separated by a vertical line,
or $(u,x] \cap(u',x'] = \varnothing$. Similarly, $B_{H_1,1/2}$ has
independent rectangular increments in the vertical direction and
$B_{1/2,1/2}$ has independent rectangular increments in arbitrary direction.
It is also clear that for $H_1 =1 $ (resp., $H_2 =1$)
$B_{H_1, H_2}$ has invariant rectangular increments in the horizontal
(resp., vertical) direction.

Let $H_i \neq1/2,1, i=1,2 $ and $\ell$ be any line passing through
the origin. Let $K =  K_{(x-1,y-1); (x,y)}, K' = K_{(0,0); (1,1)}$ be
two rectangles
whose all sides are equal to 1. Clearly, if $x$ and $y$ are large
enough, $K$ and $K'$ are separated by an
orthogonal line $\ell' \bot\ell$. From \eqref{FBsheet} and Taylor's
expansion, it easily follows that
\begin{eqnarray*}
\mathrm{E}B_{H_1,H_2}(K) B_{H_1,H_2}\bigl(K'\bigr) &
\sim&C(H_1,H_2)x^{2H_1-2} y^{2H_2-2} \qquad
\mbox{when } x, y \to \infty,
\end{eqnarray*}
with\vspace*{-3pt}
\[
C(H_1,H_2) := \prod_{i=1}^2
(2H_i) (2H_i -1) \ne0.
\]
This means that for $H_i \notin\{1/2,1\}, i=1,2$, $B_{H_1,H_2}$ has
properly dependent rectangular
increments in arbitrary direction $\ell$.

Using the terminology of Definition~\ref{iinc}, we conclude that
fractional Brownian sheet $B_{H_1,H_2}$ has:
\begin{itemize}
\item properly dependent rectangular increments if $H_i \notin\{1/2,
1\}, i=1,2$;

\item independent rectangular increments in the horizontal (vertical) direction
if $H_1 = 1/2 $ ($H_2 =1/2) $;

\item invariant rectangular increments in the horizontal (vertical) direction
if $H_1 =1 $ ($H_2 =1$);

\item independent rectangular increments if $H_1 = H_2 = 1/2$.
\end{itemize}
\end{ex}

\begin{defn} \label{Adlm}
Let $ Y = \{Y(t,s); (t,s) \in\mathbb{Z}^2
\}$ be
a stationary RF. Assume that for any $\gamma>0$ there exist a
normalization $A_n(\gamma) \to\infty$
and a non-trivial
RF  $V_\gamma= \{V_\gamma(x,y); (x,y) \in\bar\mathbb{R}^2_+\}$
such that~\eqref{Xsum} holds.

We say that $Y$ \textup{has Type I distributional LRD} (or \textup{$Y$
is a Type I RF}) if there exists
$\gamma_0 > 0$ such that
\begin{itemize}
\item RF $V_{\gamma_0} $ has properly dependent rectangular
increments, and

\item RFs $V_\gamma, \gamma\neq\gamma_0 $ do not have properly
dependent rectangular increments; in other words,
for any $\gamma\ne\gamma_0, \gamma>0$ there exists a line $\ell
(\gamma) \in\mathbb{R}^2$ such that
$V_\gamma$ has either independent or invariant increments in the
direction $\ell(\gamma)$.
\end{itemize}
Moreover, a Type I RF  $Y$ is said to have \textup{isotropic
distributional LRD} if $\gamma_0 = 1$ and
\textup{anisotropic distributional LRD} if $\gamma_0 \neq1$.
\end{defn}

\begin{rem}\label{typeI}
The above definition does not assume the
occurrence of scaling transition at $\gamma_0$, although
in all cases known to us, Type I distributional LRD property holds
simultaneously with scaling transition.
On the other hand, Remark~\ref{remON} shows that scaling transition
need not lead to Type I distributional LRD.
``Type I'' indicates that $V_{\gamma} $ has properly dependent
rectangular increments at a
\textit{single} point $\gamma= \gamma_0$. By contrast, ``Type II'' Gaussian
LRD RFs mentioned in Remark~\ref{remGaus} below have
the property that $V_{\gamma} $ have properly dependent rectangular
increments for \textit{all} $\gamma>0$.
\end{rem}

\begin{rem}\label{remGaus}
Puplinskait\.e and Surgailis \cite
{ps2014} established scaling transition and Type I distributional LRD property
for stationary Gaussian RFs with spectral density $f(x,y)= g(x,y)
(|x|^2 + |y|^{2 {\mathcal{H}}_2/{\mathcal{H}}_1} )^{-{\mathcal
{H}}_1/2}, (x,y) \in[-\pi
,\pi
]^2 $,
where ${\mathcal{H}}_i >0,  {\mathcal{H}}_1 {\mathcal{H}}_2 <
{\mathcal{H}}_1+ {\mathcal{H}}_2$ are
parameters and $g$ is a bounded positive function having nonzero limit
at the origin.
In this case, $\gamma_0 = {\mathcal{H}}_1/{\mathcal{H}}_2 $ and the
unbalanced
scaling limits $V_\pm$ agree with a fractional
Brownian sheet $B_{H_1,H_2}$ where at least one of the two parameters
$H_1, H_2 $ equals $1/2$ or 1. Moreover,
${\mathcal{H}}_1 = {\mathcal{H}}_2 $ (resp., ${\mathcal{H}}_1 \ne
{\mathcal{H}}_2$) correspond to Type
I isotropic (resp., anisotropic) distributional LRD
properties.
By contrast, ``Type II'' Gaussian RFs with spectral density of the form
$f(x,y)= g(x,y) |x|^{-2d_1}|y|^{-2 d_2},  0< d_1, d_2 < 1/2$ and a
similar function
$g$ do not exhibit scaling transition since their scaling limits
$V_\gamma$ for any $\gamma>0$
coincide with a fractional Brownian sheet $B_{d_1+ 0.5, d_2+ 0.5} $
up to a multiplicative constant; see \cite{ps2014}. \cite{lav2007,lls2014} discuss scaling limits of Gaussian LRD RFs
with general anisotropy axis.
\end{rem}

\begin{rem}\label{remON}
Scaling transition different from Type I
arises under joint temporal and contemporaneous
aggregation of \textit{independent} LRD processes in telecommunication and
economics; see \cite{miko2002,gaig2003,domb2011,pils2014} and the references therein. In these works, $\{Y(t,s);
t \in\mathbb{Z}\}, s \in\mathbb{Z}$ are independent copies of a stationary
LRD process $X = \{X(t); t \in\mathbb{Z}\}$ and the scaling limits
$V_\gamma$
of RF $Y = \{Y(t,s); (t,s)\in\mathbb{Z}^2\}$ necessarily have
independent increments in the vertical direction for any $\gamma>0$,
meaning that $Y$ cannot have
Type I distributional LRD by definition. Nevertheless, for heavy-tailed
centered ON/OFF process $X$ and some other duration based
models, the results in \cite{miko2002} imply that $Y$ exhibits a
scaling transition with some $\gamma_0 \in(0,1)$ and
markedly distinct ``supercritical'' and ``subcritical'' unbalanced scaling
limits $V_\pm$, $V_+$ being a Gaussian RF with dependent increments in
the horizontal direction
and $V_-$ having $\alpha$-stable $(1< \alpha< 2)$ distributions and
independent increments in the horizontal direction. The well-balanced
scaling limit $V_{\gamma_0}$ termed the ``intermediate process'' is
discussed in detail in \cite{gaig2006,pils2014}.
\end{rem}

\begin{proposition}\label{Vanis}
Let $ Y = \{Y(t,s); (t,s) \in\mathbb
{Z}^2\}$
be a stationary RF satisfying \eqref{Xsum} for some $\gamma>0$
and $A_n(\gamma) = L(n)n^{H}$, where $H>0$ and $L: [1,\infty) \to
\mathbb{R}_+$
is a slowly varying function.
Then the limit RF $ V_\gamma$ in \eqref{Xsum01} satisfies
the self-similarity property \eqref{OS1}. In particular, $ V_\gamma$
is OSRF corresponding
to $E := \operatorname{diag}(1, \gamma)$. Moreover, $ V_\gamma$ has stationary
rectangular increments.
\end{proposition}

\begin{pf}
Fix $\lambda>0$ and let $m:= n \lambda^{1/H}$. Then
$L(n)/L(m) \to1,  n \to\infty$ and
\begin{eqnarray*}
V_\gamma\bigl(\lambda^{1/H} x, \lambda^{\gamma/H} y\bigr)
&=&\fddlimn \frac{1}{n^{H}L(n) } \sum_{(t,s) \in K_{[x \lambda^{1/H}
n, y\lambda^{\gamma/H} n^{\gamma}]}} Y(t,s)
\\
&=&\fddlimm \frac{L(m)}{L(n)} \frac{\lambda}{m^{H} L(m)} \sum
_{(t,s) \in K_{[xm, y m^{\gamma}]}} Y(t,s) \eqfdd \lambda V_\gamma(x,y).
\end{eqnarray*}
The fact that $V_\gamma$ has stationary rectangular increments is an
easy consequence of
$Y$ being stationary.
\end{pf}

\section{Scaling transition in the aggregated 3N model}\label{sec3}

This section establishes scaling transition and Type I
anisotropic distributional LRD property, in the sense of Definitions
\ref{phase} and
\ref{Adlm} of Section~\ref{sec2}, for the aggregated 3N model ${\mathfrak X}_3$
in \eqref{aggre34}.
We shall assume that $M$ in \eqref{aggre34} is symmetric $\alpha
$-stable with
characteristic function $\mathrm{E}\mathrm{e}^{{\mathrm{i}} \theta
M(B) }
= \mathrm{e}^{ - |\theta|^\alpha\Phi(B) },  B \subset[0,1)$.
The case of general $\alpha$-stable random measure $M$ (see \eqref
{Mdef}) in \eqref{aggre34} can be discussed in a similar way.
Recall that $g_3(t,s,a)$ in \eqref{aggre34} is the Green function of
the random walk
$\{W_k \} $ on $\mathbb{Z}^2 $ with one-step transition probabilities
shown in
Figure~\ref{fig1}(a). According to Remark~\ref{alpha},
RF ${\mathfrak X}_3$ in \eqref{aggre34} with mixing distribution in
\eqref{mixdensity}
is well-defined if $1< \alpha\le2,
\beta> - (\alpha-1)/2$.

For given $\gamma>0$, introduce a RF $ V_\gamma= \{V_{3 \gamma}(x,y);
(x,y) \in\bar\mathbb{R}^2_+\}$ written as a stochastic\vspace*{-2pt} integral
%
\begin{eqnarray}
\label{3Lfield0} V_{3 \gamma}(x,y)&:=&\int_{\mathbb{R}^2 \times\mathbb{R}_+}
F_{3\gamma}(x,y; u,v,z) \mathcal{M}(\mathrm{d}u, \mathrm{d}v, \mathrm{d}z),
\end{eqnarray}
where $F_{3\gamma}(x,y; u,v,z)$ is defined as
%
\begin{eqnarray}
\label{Fgamma} &&F_{3\gamma} := %
\cases{ \displaystyle \int
_0^x \!\int_0^y
h_3(t-u, s-v, z)\, \mathrm{d}t \,\mathrm{d}s, \cr
\displaystyle \qquad \gamma = 1/2,
\cr
\displaystyle \mathbf{1}(0<v<y) \int_0^x
\mathrm{d}t \int_\mathbb{R}h_3(t-u, w, z)\,
\mathrm{d}w, \cr
\displaystyle\qquad\gamma> 1/2, 0 < \beta< \alpha-1, \cr
\displaystyle x \int_0^y h_3(-u,
s-v, z) \,\mathrm{d}s, \cr
\displaystyle \qquad \gamma> 1/2, -(\alpha -1)/2 < \beta< 0,
\cr
\displaystyle \mathbf{1}(0<u < x) \int_0^y
\, \mathrm{d}s \int_\mathbb{R}h_3(w, v-s, z) \,
\mathrm{d}w, \cr
\displaystyle \qquad \gamma< 1/2, (\alpha-1)/2 < \beta< \alpha-1, 
\cr
\displaystyle y \int_0^x
h_3(t-u, v, z) \,\mathrm{d}t, \cr
\displaystyle \qquad \gamma< 1/2, -(\alpha-1)/2 <
\beta < (\alpha-1)/2,}
\end{eqnarray}
$h_3(t,s,z) = \frac{3}{2\sqrt{\pi t}}\mathrm{e}^{-3zt - {s^2}/({4t})} \mathbf{1}
(t>0, z >0) $ as in
\eqref{hgreen}, and
$\mathcal{M}$ is an $\alpha$-stable random measure on $\mathbb{R}^2
\times
\mathbb{R}
_+$ with control measure
$\mathrm{d}\mu(u,v,z) := \phi_1 z^{\beta} \,\mathrm{d}u\, \mathrm
{d}v \,\mathrm{d}z $ and characteristic
function $\mathrm{E}\mathrm{e}^{{\mathrm{i}} \theta\mathcal{M}
(B)} = \mathrm
{e}^{-|\theta|^\alpha
\mu(B)}$,  where $\phi_1 >0, \beta> -1 $
are the asymptotic parameters in \eqref{mixdensity} and
$ B \subset\mathbb{R}^2 \times\mathbb{R}_+$ is a measurable
set with $\mu(B) < \infty$.

\begin{proposition} \label{3Nexist}
\textup{(i)} The RF $ V_{3\gamma}$
in (\ref{3Lfield0}) is well-defined for any $\gamma>0, 1< \alpha\le
2 $ and $\beta$ in~\eqref{Fgamma}.
It has $\alpha$-stable finite-dimensional distributions and stationary
rectangular increments
in the sense of \eqref{Vstatinc}.

\textup{(ii)} $V_{3\gamma}$ is OSRF: for any $\lambda>0$,
\[
\bigl\{V_{3\gamma}\bigl(\lambda x, \lambda^\gamma y\bigr); (x,y)
\in\bar\mathbb {R}^2_+\bigr\} \eqfdd\bigl\{\lambda^{H(\gamma)}
V_{3\gamma}(x,y); (x,y) \in\bar \mathbb{R} ^2_+\bigr\},
\]
where
%
\begin{equation}
\label{Hgamma} H(\gamma):= %
\cases{ \displaystyle\frac{\gamma+ \alpha- \beta}{\alpha}, &
\quad $\gamma\ge1/2, \beta >0$, \vspace*{3pt}
\cr
\displaystyle\frac{\gamma+ \alpha- 2\beta\gamma}{\alpha}, & \quad $
\gamma\ge1/2, \beta<0$,\vspace*{3pt}
\cr
\displaystyle \frac{1-\gamma+ 2\gamma(\alpha- \beta)}{\alpha}, & $
\quad \gamma< 1/2, \beta> (\alpha-1)/2$, \vspace*{3pt}
\cr
\displaystyle
\frac{\alpha\gamma+(\alpha+ 1)/2 - \beta}{\alpha}, & $\quad\gamma< 1/2, \beta< (\alpha-1)/2$.}
\end{equation}

\textup{(iii)}  RF $V_{3\gamma}$ has properly dependent rectangular
increments for $\gamma= 1/2$ and
does not have properly dependent rectangular increments for $\gamma\ne
1/2 $.

\textup{(iv)}  RFs $V_{3\gamma} = V_{3,+}  (\gamma>1/2) $ and
$V_{3\gamma}
= V_{3,-}  (\gamma<1/2) $
do not depend on $\gamma$ for $\gamma> 1/2$ and $\gamma< 1/2$.

\textup{(v)}  For $\alpha= 2$, the RFs
%
\begin{eqnarray}
V_{3,+} &\eqfdd & \kappa_{3,+} %
\cases{
B_{1- (\beta/2),1/2}, & \quad $0< \beta< 1$, \vspace*{3pt}
\cr
B_{1,(1/2) - \beta}, &
\quad $-1/2 < \beta< 0$,}
\nonumber
\\[-8pt]
\\[-8pt]
\nonumber
V_{3,-} &\eqfdd & \kappa_{3,-} %
\cases{
B_{1/2, (3/2)- \beta}, & $\quad 1/2< \beta< 1$, \vspace*{3pt}
\cr
B_{(3/4) - (\beta/2),1}, &
\quad $-1/2 < \beta< 1/2$,}
\end{eqnarray}
agree, up to some constants $\kappa_{3,\pm} = \kappa_{3,\pm} (\beta)
\neq0$, with fractional Brownian sheet $B_{H_1,H_2}$
where one of the parameters $H_1, H_2$ equals $1/2$ or 1.
\end{proposition}

\begin{rem}
Similarly, as in the case of fractional Brownian sheet
(case $\alpha=2$), the unbalanced
limit RFs $V_{3,\pm} $ have a very special dependence structure, being
either ``independent'' or ``deterministic continuations'' of random
processes with one-dimensional time:
%
\begin{eqnarray}
\mathcal{V}_{11}&:=&\bigl\{V_{3,+}(t,1); t\ge0\bigr\},
\qquad 0< \beta< \alpha-1,
\nonumber
\\
\mathcal{V}_{12}&:=&\bigl\{V_{3,+}(1,t); t\ge0\bigr\},
\qquad -(\alpha-1)/2 < \beta < 0,
\nonumber
\\[-8pt]
\label{calY}
\\[-8pt]
\mathcal{V}_{21}&:=&\bigl\{V_{3,-}(1,t); t\ge0\bigr\},
\qquad (\alpha-1)/2 < \beta < \alpha-1,
\nonumber
\\
\mathcal{V}_{22}&:=&\bigl\{V_{3,-}(t,1); t\ge0\bigr\},
\qquad -(\alpha-1)/2< \beta < (\alpha-1)/2.
\nonumber
\end{eqnarray}
The four processes $\mathcal{V}_{ij}, i,j=1,2$ in \eqref{calY} are all
symmetric $\alpha$-stable (S$\alpha$S) and self-similar with stationary
increments
(SSSI) with corresponding self-similarity parameters:
\begin{eqnarray*}
H_{11}&:=& \frac{\alpha-\beta}{\alpha},\qquad H_{12}:=
\frac{1-2\beta}{\alpha},
\\
H_{21}&:=& \frac{2(\alpha-\beta)-1}{\alpha},\qquad H_{22}:=
\frac{\alpha+1 -2\beta}{2\alpha}.
\end{eqnarray*}
These facts follow from Proposition~\ref{3Nexist}, for example, the
self-similarity property of $\mathcal{V}_{12}$ follows from
the definition of $V_{3,+}$ and Proposition~\ref{3Nexist}(ii):
$\forall\lambda>0$,
\begin{eqnarray*}
\bigl\{\mathcal{V}_{12}(\lambda t)\bigr\}&=&\bigl
\{V_{3,+}(1,\lambda t)\bigr\} = \bigl\{ V_{3,+}\bigl(
\lambda^{1/\gamma} \lambda^{-1/\gamma} 1, \lambda t\bigr)\bigr\}
\\
&\eqfdd&\lambda^{H(\gamma)/\gamma} \bigl\{V_{3,+}\bigl(
\lambda^{-1/\gamma
}1,t\bigr)\bigr\} = \lambda^{(H(\gamma)-1)/\gamma} \bigl
\{V_{3,+}(1,t)\bigr\} = \lambda^{H_{12} } \bigl\{
\mathcal{V}_{12}(t)\bigr\}.
\end{eqnarray*}
For $\alpha=2$, processes $\mathcal{V}_{ij}, i,j=1,2$ are representations
of fractional
Brownian motion and, for $1< \alpha< 2 $, they belong to the class of
S$\alpha$S SSSI processes discussed in \cite{sur1992}. Note that the
self-similarity exponents satisfy
$1/\alpha< H_{ij} <1,  i,j=1,2$ and fill in all points of the
interval $(1/\alpha,1)$ as $\beta$ vary in the corresponding intervals
in \eqref{calY}.
\end{rem}

\begin{pf*}{Proof of Proposition \protect\ref{3Nexist}}
(i)\vspace*{1pt} It suffices to show
$J_\gamma(x,y) :=
\int_{\mathbb{R}^2 \times\mathbb{R}_+} |F_{3\gamma}(x,y;
u,v,z)|^\alpha\times\mu(\mathrm{d}u,
\mathrm{d}
v, \mathrm{d}z) < \infty,  x,y > 0 $. For simplicity,
we restrict the proof to $x=y=1$, or $J_\gamma< \infty, J_\gamma:=
J_\gamma(1,1)$.

First, consider the case $\gamma= 1/2$. Write $ J_{1/2} = J' + J''$,
where $J' :=
\int_{\mathbb{R}^2 \times\mathbb{R}_+}  (\int_0^1 \!\int_0^1
h_3(t-u, s-v, z)\,
\mathrm{d}t \,\mathrm{d}
s )^\alpha
\mathbf{1}(|v| \le2) \,\mathrm{d}\mu,  J'' :=
\int_{\mathbb{R}^2 \times\mathbb{R}_+}  (\int_0^1 \int_0^1
h_3(t-u, s-v, z)\,
\mathrm{d}t \,\mathrm{d}
s )^\alpha
\mathbf{1}(|v| > 2) \,\mathrm{d}\mu$. Then\vspace*{-2pt}
\begin{eqnarray*}
J'&\le&C\int_{-\infty}^1 \mathrm{d}u
\int_0^\infty z^\beta\,\mathrm {d}z
\biggl(\int_0^1 \frac{\mathbf{1} (t >u) \,\mathrm{d}t }{\sqrt{t-u}}
\mathrm{e}^{-3z(t-u)} \biggr)^\alpha = C \biggl(\int
_{-\infty}^0 \mathrm{d}u + \cdots+ \int
_0^1 \mathrm{d}u \cdots \biggr)
\\[-2pt]
&=:& C\bigl(J'_1 + J'_2
\bigr).
\end{eqnarray*}
By Minkowski's inequality,\vspace*{-2pt}
\begin{eqnarray*}
J'_1&\le& \biggl\{\int_0^1
\mathrm{d}t \biggl(\int_{-\infty}^0
\frac
{\mathrm{d}
u}{(t-u)^{\alpha/2}} \int_0^\infty
\mathrm{e}^{-(3\alpha/2) z(t-u)} z^\beta\,\mathrm{d} z \biggr)^{1/\alpha}
\biggr\}^\alpha
\\[-2pt]
&\le& C \biggl\{\int_0^1 \mathrm{d}t \biggl(
\int_0^{\infty} \frac
{\mathrm{d}u}{(t+u)^{1+
\beta+ (\alpha/2)}}
\biggr)^{1/\alpha} \biggr\}^\alpha= C \biggl\{ \int
_0^1 \frac{\mathrm{d}t}{t^{(1/2)+ (\beta/\alpha)}} \biggr\}^{1/\alpha}
< \infty
\end{eqnarray*}
since $(1/2)+ (\beta/\alpha)< 1 $
due to $\beta< \alpha-1, \alpha\le2$. We also have\vspace*{-2pt}
\[
J'_{2}\le
C\int_0^\infty z^\beta\,\mathrm{d}z  \Biggl\{ \int_0^1 \mathrm{e}^{-(3\alpha/2)zx}
\,\mathrm{d}x  \Biggr\}
^\alpha =
C \int_0^\infty z^{\beta-\alpha}  (1 - \mathrm{e}^{-z}
)^\alpha\,\mathrm{d}z  <
 \infty
 \]
since $\alpha> 1 + \beta$. On the other hand, since $(s-v)^2 \ge v^2/4
$ for $|s| < 1,  |v| > 2 $, so using
Minkowski's inequality we obtain
\begin{eqnarray*}
J'' &\le& \biggl\{\int_0^1
\mathrm{d}t \biggl(\int_{-\infty}^t
\frac
{\mathrm{d}
u}{(t-u)^{\alpha
/2}} \int_{|v|>2} \mathrm{e}^{- v^2/4(t-u)}
\,\mathrm{d}v \int_0^\infty\mathrm{e}^{-(3\alpha/2) z(t-u)}
z^\beta\,\mathrm{d}z \biggr)^{1/\alpha} \biggr\}^\alpha
\\[-2pt]
&\le&C\int_0^{\infty} \frac{\mathrm{d}x}{x^{1+ \beta+ (\alpha
/2)}} \int
_{|v|>2} \mathrm{e}^{- v^2/(4x)} \,\mathrm{d}v,
\end{eqnarray*}
where the last integral is easily seen to be finite. This proves
$J_{1/2}< \infty$.

Next, consider $J_{\gamma} $ for $\gamma> 1/2, 0< \beta< \alpha-1$.
Using $h_\star(u,z) := \int_\mathbb{R}h_3(u,v,z) \,\mathrm{d}v = 12
\mathrm
{e}^{-3uz} \mathbf{1}(u>0)$, similarly as above we\vspace*{-2pt} obtain
\begin{eqnarray*}
J_\gamma&\le&C\int_{-\infty}^1 \mathrm{d}u
\int_0^\infty z^\beta \mathrm{d}z
\biggl(\int_{u \vee0}^1 \mathrm{e}^{-3z(t-u)}
\,\mathrm{d}t \biggr)^\alpha = C \biggl\{ \int_{-\infty}^0
\mathrm{d}u +\cdots+ \int_0^1 \mathrm{d} u
\cdots \biggr\}
\\
&=:& C \{J_{\gamma1} + J_{\gamma2} \},
\end{eqnarray*}
where
$J_{\gamma1}\le C \{\int_0^1 \mathrm{d}t  (\int_0^{\infty}
(t+u)^{-1-\beta} \mathrm{d}u  )^{1/\alpha}  \}^\alpha
\le C \{\int_0^1 t^{-\beta/\alpha} \,\mathrm{d}t  \}
^{1/\alpha} <
\infty
$\vspace*{-2pt} and
\[
J_{\gamma2} \le C \int_0^1 \mathrm{d}u
\int_0^\infty z^\beta \,\mathrm{d}z
\biggl(\int_u^1 \mathrm{e}^{-3z (t-u)}
\,\mathrm{d}t \biggr)^\alpha \le C \int_0^\infty
z^\beta\,\mathrm{d}z \bigl(\bigl(1 - \mathrm {e}^{-z}\bigr)/z
\bigr)^\alpha < \infty
\]
because of $\beta- \alpha< - 1$. This proves $J_\gamma< \infty$ for
$ \gamma> 1/2, 0< \beta< \alpha-1$.

Next, let $\gamma> 1/2, -(\alpha-1)/2 < \beta< 0$. We have
\begin{eqnarray*}
J_\gamma&\le&C\int_{\mathbb{R}_+ \times\mathbb{R}\times\mathbb
{R}_+} \,\mathrm{d}\mu \biggl(\int
_0^1 h_3(u,s-v,z) \,\mathrm{d}s
\biggr)^\alpha
\\
&\le&C\int_0^\infty u^{-\alpha/2} \mathrm{d}u
\int_{\mathbb{R}} \mathrm{d}v \int_0^\infty
\mathrm{e} ^{-zu} z^\beta\,\mathrm{d}z \biggl(\int
_0^1 \mathrm{e}^{-(s-v)^2/u} \,\mathrm{d}s
\biggr)^\alpha
\\
&=&C\int_0^\infty u^{-(1 + \beta+ {\alpha}/{2})} \,\mathrm{d}u
\biggl\{ \int_{|v|\le2} \mathrm{d}v + \int_{|v|>2}
\mathrm{d}v \biggr\} \biggl(\int_0^1
\mathrm{e} ^{-(s-v)^2/u}\, \mathrm{d}s \biggr)^\alpha =: C
\{J_{\gamma1} + J_{\gamma
2} \}.
\end{eqnarray*}
Here,
\begin{eqnarray*}
J_{\gamma1}&\le& C\int_0^\infty
u^{-(1 + \beta+ {\alpha}/{2})}\, \mathrm{d}u \biggl(\int_0^1
\mathrm{e} ^{-s^2/u} \,\mathrm{d}s \biggr)^\alpha
\\
&\le& C \biggl(\int_0^1 u^{-(1 + \beta)}
\,\mathrm{d}u + \int_1^\infty u^{-(1 +
\beta+ {\alpha}/{2})}
\,\mathrm{d}u \biggr) < \infty
\end{eqnarray*}
since $\beta<0, \beta> - \alpha/2$, while
\begin{eqnarray*}
J_{\gamma2}&\le& C\int_0^\infty
u^{-(1 + \beta+ {\alpha}/{2})} \,\mathrm{d}u \int_1^\infty
\mathrm{e}^{-v^2/u} \,\mathrm{d}v \\
&\le & C \int_0^\infty
u^{-({1}/{2} + \beta+ {\alpha}/{2})} \,\mathrm{d}u \int_{u^{1/2}}^\infty
\mathrm{e}^{-z^2} \,\mathrm{d}z < \infty
\end{eqnarray*}
as $\frac{1}{2} + \beta+ \frac{\alpha}{2} > 1$ and $\int_{1/u^{1/2}}^\infty\mathrm{e}^{-z^2} \,\mathrm{d}z $ decays\vspace*{1pt}
exponentially when
$u \to0$.
This proves $J_\gamma< \infty$ for $ \gamma> 1/2, -(\alpha-1)/2 <
\beta< 0$.

Consider the case $0< \gamma< 1/2, (\alpha-1)/2 < \beta< \alpha-1 $.
Then using $
\int_\mathbb{R}h_3(w, v , z) \,\mathrm{d}w = \frac{\sqrt{3} }{2\sqrt
{z}} \mathrm
{e}^{-\sqrt{3z}
|v|} $ we obtain
\begin{eqnarray*}
J_\gamma&\le&C\int_{\mathbb{R}} \,\mathrm{d}v \int
_0^\infty z^\beta\, \mathrm{d}z \biggl(
\int_{0}^1 z^{-1/2}
\mathrm{e}^{-\sqrt{z}|s-v|} \,\mathrm{d}s \biggr)^\alpha = C \biggl\{ \int
_{|v|
\le2} \mathrm{d}v +\cdots+ \int_{|v|>2}
\mathrm{d}v \cdots \biggr\}
\\
&=:& C \{J_{\gamma1} + J_{\gamma2} \},
\end{eqnarray*}
where $J_{\gamma1} \le C \int_0^\infty z^{\beta- (\alpha/2)}
\,\mathrm{d}z
(\int_{0}^1 \mathrm{e}^{-z|s|}\, \mathrm{d}s )^\alpha
\le C \int_0^\infty z^{\beta- \alpha} (1- \mathrm{e}^{-\sqrt
{z}})^\alpha\,\mathrm{d}z
< \infty$ for $0 < \beta< \alpha-1 $ and
\begin{eqnarray*}
J_{\gamma2}&\le&C\int_1^\infty\mathrm{d}v
\int_0^\infty z^{\beta
- \alpha/2}
\mathrm{e}^{-\sqrt{z} v}\, \mathrm{d}z = C \int_1^\infty
v^{\alpha- 2 - 2\beta}\, \mathrm{d}v < \infty
\end{eqnarray*}
since $2 + 2\beta- \alpha> 1 $ for $\beta> (\alpha-1)/2$.

Finally, let $0< \gamma< 1/2, -(\alpha-1)/2 < \beta< (\alpha-1)/2$. Then
$J_\gamma= C \int_{-\infty}^1 \mathrm{d}u \int_{\mathbb{R}}
\mathrm{d}v \int_0^\infty
z^\beta\,\mathrm{d}z  \times (\int_0^1 h_3(t-u,v,z) \,\mathrm{d}t
)^\alpha
= C  \{ \int_{-\infty}^0 \mathrm{d}u +\cdots+ \int_0^1 \mathrm
{d}u \cdots
\}
=: C \{ J_{\gamma1} + J_{\gamma2}\}$. By Minkowski's
inequality,
\begin{eqnarray*}
J^{1/\alpha}_{\gamma1}&\le&C \int_0^1
\mathrm{d}t \biggl\{\int_0^\infty\mathrm{d}u \int
_0^\infty\mathrm{d}v \int_0^\infty
h^\alpha_3(t+u, v, z) z^\beta\,\mathrm{d}z \biggr
\}^{1/\alpha}
\\
&=&C \int_0^1 \mathrm{d}t \biggl\{\int
_0^\infty\frac{\mathrm{d}u}{
(t+u)^{1 +
\beta+
(\alpha-1)/2}} \biggr
\}^{1/\alpha} = \int_0^1 \mathrm{d}t \biggl
\{ \frac{1}{t^{\beta+ (\alpha
-1)/2}} \biggr\} ^{1/\alpha} < \infty
\end{eqnarray*}
and, similarly,
\begin{eqnarray*}
J^{1/\alpha}_{\gamma2}&\le&C \int_0^1
\mathrm{d}t \biggl\{\int_0^\infty\mathrm{d}v \int
_0^\infty h^\alpha_3(t, v,
z) z^\beta\,\mathrm{d}z \biggr\}^{1/\alpha} = C\int
_0^1 \mathrm{d}t \biggl\{ \frac{1}{t^{\beta+ (\alpha
+1)/2}}
\biggr\} ^{1/\alpha} < \infty
\end{eqnarray*}
since $|\beta| < (\alpha-1)/2$. This proves $J_\gamma< \infty$, or
the existence of $V_{3\gamma} $,
for all choices of $\alpha, \beta, \gamma$ in \eqref{Fgamma}.
The fact that linear combinations
of integrals in \eqref{3Lfield0} are $\alpha$-stable is well known
(\cite{sam1994}). Stationarity of increments
of \eqref{3Lfield0} is an easy consequence of the integrand \eqref
{Fgamma} and the control measure $\mu$.
This proves part (i).

(ii) The OSRF property is immediate from the scaling properties
$h_3(\lambda u, \sqrt{\lambda} v, \lambda^{-1} z)
= \lambda^{-1/2} h_3(u,v,z) $ of the kernel $h_3$ in \eqref{hgreen} and
$\{\mathcal{M}(\mathrm{d}\lambda u, \mathrm{d}\lambda^\gamma v,
\mathrm{d}\lambda^{-1}
z)\}
\eqfdd\{\lambda^{({\gamma-\beta})/{\alpha}} \mathcal
{M}(\mathrm{d}u,\break
\mathrm{d}v,
\mathrm{d}z)\} $ of the stable random measure $\mathcal{M}$,
the last property being a consequence of the scaling property
of $\mu(\mathrm{d}\lambda u, \mathrm{d}\lambda^\gamma v, \mathrm
{d}\lambda^{-1} z) =
\lambda
^{\gamma-\beta} \mu(\mathrm{d}u, \mathrm{d}v, \mathrm{d}z)$
of the control measure $\mu$.

(iii) Let $\gamma= \gamma_0 := 1/2 $.
Consider\vspace*{1pt} arbitrary rectangles $K_i = K_{(\xi_i,\eta_i); (x_i,y_i)}
\subset\mathbb{R}^2_+, i=1,2$,
and write\vspace*{1pt} $\int= \int_{\mathbb{R}^2 \times\mathbb{R}_+}$. Then
$V_{3 \gamma_0}(K_i) = \int G_{K_i}(u,v,z) \,\mathrm{d}\mathcal{M}$,
where $G_{K_i}(u,v,z) := \int_{K_i} h_3(t-u, s-v, z) \,\mathrm{d}t
\,\mathrm{d}s$.
Note $G_{K_i} \ge0$ and
$G_{K_i}(u,v,z) >0$ for any $u < x_i$ implying $\operatorname{supp}(G_{K_1})
\cap\operatorname{supp}(G_{K_2}) \neq\varnothing$.
Hence, and from (\cite{sam1994}, Theorem~3.5.3, page~128)
it follows that the
increments $V_{3 \gamma_0}(K_i), i=1,2 $ on arbitrary nonempty
rectangles $K_1, K_2 $ are dependent. It is also
easy to show that $V_{3\gamma_0}$ does not have invariant rectangular
increments
in any direction. This proves (iii) for $\gamma= 1/2$.

Next,\vspace*{1.5pt} let $\gamma> 1/2, 0< \beta< \alpha-1$. Similarly\vspace*{1pt} as above,
for any rectangle $K = K_{(\xi,\eta); (x,y)} \subset\mathbb{R}^2_+$, we have
$V_{3\gamma} (K) = \int G_{K,\gamma}(u,v,z) \,\mathrm{d}\mathcal{M}$,
where\vspace*{1.5pt} $G_{K,\gamma}(u,v,z) := \mathbf{1}(\eta< v \le y) \int_{\xi
}^\eta
h_{3\gamma}(t-u,z) \,\mathrm{d}t $. Clearly,
if $K_i, i=1,2 $ are any two rectangles separated by a horizontal line, then
$\operatorname{supp}(G_{K_1,\gamma}) \cap\operatorname{supp}(G_{K_2}) = \varnothing$, implying
independence of $V_{3\gamma} (K_1) $ and $V_{3\gamma}(K_2)$. Thus,
$V_{3\gamma}$ for $0< \beta< \alpha-1$ has
independent increments in the vertical direction.
The fact that
$V_{3\gamma} $ for $\gamma>1/2, -(\alpha-1)/2 < \beta<0$ has
invariant increments in the horizontal direction
is obvious from \eqref{3Lfield0} and \eqref{Fgamma}. The properties of
$V_{3\gamma} $ in the
case $0< \gamma< 1/2 $ are completely analogous.

(iv) Follows from \eqref{3Lfield0} and \eqref{Fgamma}.

(v) Since $V_{3, \pm}$ for $\alpha=2$ are zero mean Gaussian
RFs, it suffices to show that their covariances agree
with that of fractional Brownian sheet in \eqref{covB}. This can be
easily verified by using self-similarity and stationarity
of increments properties stated in (i) and (ii), as follows.

Let $0< \beta< 1$ and $\rho_+(x,x') := \mathrm{E}V_{3,+}(x,1) V_{3,+}(x',1),
x,x' \ge0$. By \eqref{3Lfield0} and \eqref{Fgamma},
$\mathrm{E}V_{3,+}(x,y) V_{3,+}(x',y') = (y \wedge y') \rho_+ (x,x'),
(x,y), (x',y') \in\mathbb{R}^2_+$. According to (ii),
for any $\lambda>0$
%
\begin{eqnarray}
\rho_+ \bigl(\lambda x, \lambda x'\bigr)&=&\mathrm{E}V_{3,+}(
\lambda x,1) V_{3,+}\bigl(\lambda x',1\bigr) =
\lambda^{2H(\gamma)} \mathrm{E}V_{3,+}\bigl(x,\lambda^{-\gamma}
\bigr) V_{3,+}\bigl(x',\lambda ^{-\gamma}\bigr)
\nonumber
\\[-8pt]
\label{rhoscale}
\\[-8pt]
\nonumber
&=&\lambda^{2H(\gamma)-\gamma} \mathrm{E}V_{3,+}(x,1)
V_{3,+}\bigl(x',1\bigr) = \lambda^{2H_+} \rho
\bigl(x,x'\bigr),
\end{eqnarray}
where $H_+ := H(\gamma) - (\gamma/2) = 1 - (\beta/2)$; see \eqref
{Hgamma}. The stationarity of rectangular increments
property of RF $V_{3,+}$ implies that the process
$\{ V_{3,+}(x,1), x \ge0 \}$ has stationary increments. Together with
the scaling property in \eqref{rhoscale},
this implies that $\rho_+ (x,x') = (\kappa_+^2/2) (x^{2H_+} + x'^{2H_+}
- |x-x'|^{2H_+}),  x,x' \ge0$, or
$\mathrm{E}V_{3,+}(x,y) V_{3,+}(x',y') = \kappa_+^2 \mathrm{E}B_{1-
(\beta/2),
1/2}(x,y) B_{1- (\beta/2), 1/2}(x',y') $, see
\eqref{covB}. The remaining relations in (v) are analogous.
Proposition~\ref{3Nexist} is proved.
\end{pf*}

The main result of this section is Theorem~\ref{3Ntheo}. Its proof is
based on the asymptotics of the Green
function $g_3$ in Lemma~\ref{lemma3N}, below. The proof of Lemma~\ref
{lemma3N} can be found at \surl{http://arxiv.org/abs/1303.2209v3}.

\begin{lem} \label{lemma3N}
For any $ (t,s,z) \in(0,\infty) \times\mathbb{R}\times(0,\infty)
$ the
point-wise convergence in \eqref{green03} holds.
This convergence is uniform on any relatively compact set $\{\epsilon<
t < 1/\epsilon,  \epsilon< |s| < 1/\epsilon,
\epsilon< z < 1/\epsilon\} \subset(0,\infty) \times\mathbb
{R}\times
(0,\infty),  \epsilon>0$.

Moreover, there exist constants $C, c >0$ such that for all
sufficiently large $\lambda$ and
any $(t,s,z),  t>0,  s\in\mathbb{R},  0 < z < \lambda$ the following
inequality holds:
%
\begin{equation}
\label{limg3bdd} \sqrt{\lambda} g_3 \biggl([\lambda t], [\sqrt{
\lambda} s], 1- \frac{z}{\lambda} \biggr) < C \bigl(\bar h_3(t,s, z)
+ \sqrt {\lambda} \mathrm{e}^{- zt - c (\lambda t)^{1/3} - c(\sqrt{\lambda}
|s|)^{1/2}} \bigr),
\end{equation}
where $\bar h_3(t,s,z) := \frac{1}{\sqrt{t}} \mathrm{e}^{-zt -
{s^2}/({16t})},  (t,s,z) \in(0,\infty) \times\mathbb{R}\times(0,\infty)$.
\end{lem}

\begin{theorem} \label{3Ntheo}
Assume that the mixing density $\phi$ is bounded on any interval
$[0,1-\epsilon), \epsilon>0$ and satisfies \eqref{mixdensity}, where
%
\begin{equation}
\label{bcond} -(\alpha-1)/2 < \beta< \alpha- 1,\qquad 1< \alpha\le2, \beta
\neq0, \beta\neq(\alpha-1)/2.
\end{equation}
Let $\mathfrak{X}_3 $ be the aggregated RF in \eqref{aggre34}.
Then for any $\gamma>0 $
%
\begin{eqnarray}
\label{3Nlim0} && n^{-H(\gamma)} \sum_{t=1}^{[nx]}
\sum_{s=1}^{[n^\gamma y]} \mathfrak
{X}_3(t,s) \limfdd V_{3\gamma}(x,y),\qquad x, y >0, n \to
\infty,
\end{eqnarray}
where $H(\gamma)$ and $V_{3\gamma}$ are given in \eqref{Hgamma} and
\eqref{3Lfield0}, respectively. As a consequence, the RF $\mathfrak
{X}_3$ exhibits scaling transition at $\gamma_0 = 1/2$ and
enjoys
Type I anisotropic distributional LRD with $\gamma_0 = 1/2$ in the
sense of Definition~\ref{Adlm}.
\end{theorem}

\begin{rem}
As it follows from the proof of Theorem~\ref{3Ntheo},
for $\gamma= 1/2 $ the limit in
\eqref{3Nlim0} exists also when $\beta= 0$ or $\beta= (\alpha-1)/2$
and is given in \eqref{3Lfield0} as in
the remaining cases. On the other hand, the existence of the scaling
limit \eqref{3Nlim0}
in the cases $\gamma>1/2, \beta=0$ and $0<\gamma< 1/2 $ and $\beta=
(\alpha-1)/2$ is an open and delicate
question. Note a sharp transition in the dependence structure of the
limit fields $V_{3,+}$ and $V_{3,-}$
in the vicinity of $\beta= 0$ and $\beta= (\alpha-1)/2 $, respectively,
changing abruptly from
independent rectangular increments in one direction to invariant
(completely dependent) rectangular increments
in the perpendicular direction. For $\alpha= 2$, the above transition
may be related to the fact that
the covariance functions
of the ``vertical'' and ``horizontal sectional processes'' $\{\mathfrak
{X}_3(0,s); s \in\mathbb{Z}\} $ and $\{\mathfrak{X}_3(t,0); t \in
\mathbb{Z}\} $
change their summability properties at respective points $\beta=0$ and
$\beta= 1/2$;
see Proposition~\ref{3Ncov} below.
\end{rem}

Let $\alpha=2$ and $ r_3 (t,s)= \mathrm{E}\mathfrak{X}_3(t,s)
\mathfrak
{X}_3(0,0)$ be the covariance function of the aggregated Gaussian RF in
(\ref{aggre34}).
The proof of Proposition~\ref{3Ncov} using Lemma~\ref{lemma3N}
can be found in the arXiv version \surl{http://arxiv.org/abs/1303.2209v3}.

\begin{proposition} \label{3Ncov}
Assume $\alpha=2$ and the conditions
of Theorem~\ref{3Ntheo}.
Then for any $(t,s) \in\mathbb{R}^2_0$
\begin{eqnarray*}
\label{3Ncovasy} \lim_{\lambda\to\infty} \lambda^{\beta+ 1/2}
r_3\bigl([\lambda t], [\sqrt {\lambda} s]\bigr)&=& %
\cases{ C_3 |s|^{-2\beta-1} \gamma\bigl(\beta+ 1/2,
s^2/4|t|\bigr), & \quad $t\ne0, s \ne0$,\vspace*{3pt}
\cr
C_3 |s|^{-2\beta-1} \Gamma(\beta+ 1/2), &\quad $t=0$,
\vspace*{3pt}
\cr
C_4 |t|^{-\beta- 1/2}, & \quad $s=0$,}
\end{eqnarray*}
where $\gamma(\alpha, x) := \int_0^x y^{\alpha-1} \mathrm{e}^{-y}
\,\mathrm{d}y $ is
incomplete gamma function and
$C_3 := \pi^{-{1}/{2}} 2^{2\beta-1}3^{1-\beta}\sigma^2\times\phi
_1\Gamma
(\beta+1)$,  $C_4 := 4^{-{1}/{2}-\beta}C_3$.
\end{proposition}

\begin{pf*}{Proof of Theorem \protect\ref{3Ntheo}}
Write $S_{n\gamma}(x,y) $ for the left-hand side of \eqref{3Nlim0}. It suffices
to prove the convergence
of characteristic functions:
%
\begin{equation}
\label{fidi} \mathrm{E}\mathrm{e}^{\mathrm{i}\sum_{j=1}^p \theta_j S_{n\gamma
}(x_j,y_j) } \to \mathrm{E}\mathrm{e}
^{\mathrm{i}\sum_{j=1}^p \theta_j V_{3\gamma}(x_j,y_j) },\qquad n \to \infty,
\end{equation}
for any $p\in\mathbb{N}_+, \theta_j \in\mathbb{R}, (x_j,y_j) \in
\mathbb{R}^2_+, j=1,
\ldots, p$.
We have
%
\begin{equation}
\mathrm{E}\mathrm{e}^{\mathrm{i}\sum_{j=1}^p \theta_j S_{n\gamma
}(x_j,y_j) } = \mathrm{e}^{-
J_{n\gamma}},\qquad
\mathrm{E}\mathrm{e}^{\mathrm{i}\sum_{j=1}^p \theta_j V_{3\gamma
}(x_j,y_j) } = \mathrm{e}^{-
J_\gamma},
\end{equation}
where
%
\begin{eqnarray}
J_\gamma&:=&\int_{\mathbb{R}^2 \times\mathbb{R}_+}
\bigl|G_\gamma (u,v,z)\bigr|^\alpha\,\mathrm{d} \mu,\qquad
G_\gamma(u,v,z) := \sum_{j=1}^p
\theta_j F_{3\gamma}(x_j,y_j;
u,v,z),\hspace*{6pt}
\nonumber
\\[-8pt]
\label{Jgamma}
\\[-8pt]
\nonumber
J_{n\gamma}&:=&n^{-H(\gamma) \alpha} \sum_{(u,v) \in\mathbb{Z}^2}
\mathrm{E} \Biggl|\sum_{j=1}^p \theta
_j \sum_{1\le t
\le
[nx_j], 1\le s \le[n^\gamma y_j]} g_3(t-u,
s-v, a) \Biggr|^\alpha.
\nonumber
\end{eqnarray}
Thus, \eqref{fidi} 
follows from
%
\begin{equation}
\label{fidi1} \lim_{n\to\infty} J_{n\gamma} =
J_\gamma.
\end{equation}
To prove \eqref{fidi1}, we write $J_{n\gamma} $ as an integral
%
\begin{equation}
\label{Jngamma} J_{n\gamma} = \int_{\mathbb{R}^2 \times\mathbb{R}_+}
\bigl|G_{n\gamma
}(u,v,z)\bigr|^\alpha \chi _n(z) \mu(
\mathrm{d}u, \mathrm{d}v, \mathrm{d}z),
\end{equation}
where the functions $\chi_n$ satisfying
$\chi_n(z) \to1  (n \to\infty)$ uniformly in $z >0$ will be
specified later, and where $G_{n\gamma}: \mathbb{R}^2 \times\mathbb
{R}_+ \to\mathbb{R}$ are
some functions
which approach $G_\gamma$ in \eqref{Jgamma} in the following sense. Let
$W_\epsilon:= \{ (u,v,z) \in\mathbb{R}^2 \times\mathbb{R}_+: |u|+
|v| < 1/\epsilon,
\epsilon<z< 1/\epsilon\},
W_\epsilon^c := (\mathbb{R}^2 \times\mathbb{R}_+) \setminus
W_\epsilon,
\epsilon
>0$. We will prove that
%
\begin{equation}
\label{gamma0} \lim_{n \to\infty} \int_{W_\epsilon}
\bigl|G_{n\gamma}(u,v,z) - G_\gamma (u,v,z)\bigr|^\alpha\,\mathrm{d}
\mu = 0 \qquad \forall \epsilon>0,
\end{equation}
and\vspace*{-3pt}
%
\begin{equation}
\label{gamma2} \lim_{\epsilon\to0} \limsup_{n \to\infty}
\int_{W_\epsilon^c} \bigl|G_{n\gamma}(u,v,z)\bigr|^\alpha
\,\mathrm{d}\mu = 0.
\end{equation}
Since $\mu(W_\epsilon) < \infty$, \eqref{gamma0} follows from the
uniform convergence
%
\begin{equation}
\label{gamma1} \lim_{n \to\infty} \sup_{(u,v,z) \in W_\epsilon}
\bigl|G_{n\gamma
}(u,v,z) - G_\gamma(u,v,z)\bigr| = 0 \qquad \forall
\epsilon>0.
\end{equation}
Clearly, \eqref{gamma0} and \eqref{gamma2} together with \eqref
{Jngamma} and the above mentioned property of $\chi_n$ imply~\eqref{fidi1}.

The subsequent proof of \eqref{gamma0} and \eqref{gamma2} is split into
several cases depending on values $\gamma$ and $\beta$.

\textit{Case}  $\gamma= \gamma_0 = 1/2$.  In this case, \eqref
{Jngamma} holds with\vspace*{-3pt}
%
\begin{eqnarray}
&& G_{n\gamma_0}(u,v,z)
\nonumber
\\[-8pt]
\label{Gamma0}
\\[-8pt]
\nonumber
&&\quad:=\sum_{j=1}^p
\theta_j \int_0^{\lfloor
nx_j\rfloor /n} \!\int
_0^{\lfloor \sqrt n y_j\rfloor /\sqrt n} \sqrt{n} g_3 \biggl(\lceil
nt\rceil -\lceil nu\rceil , \lceil \sqrt {n}s\rceil -\lceil \sqrt{n}v\rceil , 1 -
\frac{z}{n} \biggr) \,\mathrm{d}t \,\mathrm{d}s
\end{eqnarray}
and $\chi_n(z) := (z/n)^{-\beta} (\phi(1 - z/n)/\phi_1) \mathbf{1} (0<
z < n)
\to1 $ boundedly on $\mathbb{R}_+$ as $n \to\infty$ according to condition
\eqref{mixdensity}.
To show \eqref{gamma1}, for given $\epsilon_1 >0$ split
$G_{n\gamma_0}(u,v,z)- G_{\gamma_0}(u,v,z)  =   \sum_{i=1}^3
\Gamma
_{ni}(u,v,z) $,
where, for $0< z < n$,
\begin{eqnarray*}
\Gamma_{n1}(u,v,z)&:=&\sum_{j=1}^p
\theta_j \int_0^{{\lfloor
nx_j\rfloor }/{n}} \!\int
_0^{{\lfloor \sqrt n y_j\rfloor }/{\sqrt n}} \biggl\{ \sqrt{n} g_3
\biggl(\lceil nt\rceil - \lceil nu \rceil , \lceil \sqrt{n}s\rceil -\lceil \sqrt
{n}v\rceil , 1 - \frac{z}{n} \biggr)
\\[-2pt]
&&{}- h_3(t-u, s-v, z) \biggr\} \mathbf{1}\bigl((t,s)\in
D_j(\epsilon _1)\bigr)\,\mathrm{d}t \,\mathrm{d}s,
\\[-2pt]
\Gamma_{n2}(u,v,z)&:=&\sum_{j=1}^p
\theta_j \int_0^{\lfloor nx_j\rfloor/n} \!\int
_0^{{\lfloor \sqrt n y_j\rfloor }/{\sqrt n}} \sqrt{n} g_3 \biggl(\lceil
nt\rceil -\lceil nu \rceil , \lceil \sqrt{n}s \rceil - \lceil \sqrt {n}v\rceil ,
1 - \frac{z}{n} \biggr)
\\[-2pt]
&&{}\times \mathbf{1}\bigl((t,s)\notin D_j(\epsilon_1)
\bigr) \,\mathrm{d}t \,\mathrm{d}s,
\\[-2pt]
\Gamma_{n3}(u,v,z)&:=&-\sum_{j=1}^p
\theta_j \int_0^{{\lfloor nx_j\rfloor }/{n}} \!\int
_0^{{\lfloor \sqrt n y_j\rfloor }/{\sqrt n}} h_3(t-u, s-v, z)\mathbf{1}
\bigl((t,s)\notin D_j(\epsilon_1)\bigr) \,\mathrm{d}t
\,\mathrm{d}s,
\end{eqnarray*}
and where the sets $D_j(\epsilon_1), j=1,\ldots, p$ (depending on $u,
v$) are defined by
\[
D_j(\epsilon_1) := \bigl\{ (t,s) \in(0,x_j]
\times(0,y_j]: t-u > \epsilon _1, |s-v| >
\epsilon_1\bigr\}.
\]
Relation \eqref{gamma1} follows from\vspace*{-2pt}
%
\begin{eqnarray}
\label{gamma11} \lim_{n \to\infty} \sup_{(u,v,z) \in W_\epsilon} \bigl
\llvert \Gamma _{n1}(u,v,z)\bigr\rrvert &=&0,
\\[-2pt]
\label{gamma12} \lim_{\epsilon_1 \to0} \limsup_{n\to\infty}
\sup_{(u,v,z) \in
W_\epsilon} \bigl\llvert \Gamma_{ni}(u,v,z)\bigr
\rrvert &=&0,\qquad i=2,3.
\end{eqnarray}
Here, \eqref{gamma11} follows from Lemma~\ref{lemma3N}.
Next,
$|\Gamma_{n3}(u,v,z)| \le C\int_0^{\epsilon_1} t^{-1/2}\, \mathrm{d}t
+ C\int_{\epsilon_1}^1 t^{-1/2}\, \mathrm{d}t \times\int_{|s| < \epsilon_1}
\,\mathrm{d}s = O(\sqrt{\epsilon_1})$, implying \eqref{gamma12} for $i=3$.
Similarly, using
(\ref{limg3bdd}) we obtain
$|\Gamma_{n2}(u,v,z)| \le C \sqrt{\epsilon_1} +
C\sqrt{n} \int_0^{1} \mathrm{e}^{- c (n t)^{1/3}}\,\mathrm{d}t
\le C \sqrt{\epsilon_1} + C/\sqrt{n}$.
This proves \eqref{gamma12} for $i=2$,
and hence \eqref{gamma1}, too.

Consider\vspace*{1pt} \eqref{gamma2}. W.l.g., we can assume $p= 1,  \theta_1 = x_1
= y_1 = 1$.
With \eqref{Gamma0} and \eqref{limg3bdd} in mind, we have
$0\le G_{n\gamma_0}(u,v,z) \le C(\bar G(u,v,z) + \widetilde
G_{n}(u,v,z))$, where
\begin{eqnarray*}
\bar G(u,v,z)&:=&\int_0^{1} \!\int
_0^{1} \bar h_3(t-u, s-v, z)\,
\mathrm{d}t \,\mathrm{d}s,
\\
\widetilde G_{n}(u,v,z)&:=&\sqrt{n} \mathbf{1}(0< z < n) \int
_0^{1} \!\int_0^{1}
\mathrm{e}^{- z(t-u) - c (n (t-u))^{1/3} - c(\sqrt{n} |s-v|)^{1/2}} \mathbf{1}(t>u) \,\mathrm{d}t \,\mathrm{d}s,
\end{eqnarray*}
where $c>0$ is the same as in \eqref{limg3bdd}. Relation \eqref{gamma2}
with $G_{n\gamma}$ replaced by
$\bar G$ follows from $\bar G \in L^\alpha(\mu)$ (see Proposition~\ref{exist}, proof of (i)), since $\bar h_3(t, s, z)$
and $h_3(t, s, z)$ differ only in constants. Thus, \eqref{gamma2}
follows from
%
\begin{eqnarray}
\label{3NJ} \widetilde J_{n}&:=&\int_{\mathbb{R}^2 \times\mathbb{R}_+}
\bigl(\widetilde G_{n}(u,v,z)\bigr)^\alpha\,\mathrm{d}\mu =
o(1), \qquad n \to\infty.
\end{eqnarray}
Split $\widetilde J_{n} = \sum_{i=1}^3 I_{ni}$, where
\begin{eqnarray*}
I_{n1} &:=& \int_{(-\infty, 0] \times\mathbb{R}\times\mathbb{R}_+} (\widetilde
G_n)^\alpha\,\mathrm{d}\mu,\\
 I_{n2} &:=& \int
_{(0,1] \times[-2,2] \times\mathbb{R}_+} (\widetilde G_n)^\alpha
\,\mathrm{d}\mu,
\\
I_{n3} &:=& \int_{(0,1] \times[-2,2]^c \times\mathbb{R}_+} (\widetilde
G_n)^\alpha\,\mathrm{d}\mu,
\end{eqnarray*}
$[-2,2]^c := \mathbb{R}\setminus[-2,2]$.
Using the fact that $\int_{\mathbb{R}} \mathrm{e}^{-c
n^{1/4}|s-v|^{1/2}} \,\mathrm{d}
v =
C/\sqrt{n} $ and Minkowski's inequality,
\begin{eqnarray*}
I_{n1}&\le& Cn^{\alpha/2} \biggl\{\int_{(0,1]^2}
\mathrm{d}t \,\mathrm{d}s \biggl( \int_{\mathbb{R}_+ \times\mathbb{R}\times\mathbb{R}_+ }
\mathrm{e}^{- \alpha z(t+u) - c\alpha(n (t+u))^{1/3} - c\alpha(\sqrt{n}
|s-v|)^{1/2}}
\\
&& {}\times z^\beta\,\mathrm{d}u \,\mathrm{d}v \,\mathrm{d}z
\biggr)^{1/\alpha} \biggr\} ^{\alpha}
\\
&\le& Cn^{({\alpha-1})/{2}} \biggl\{\int_0^1
\mathrm{d}t \biggl(\int_0^\infty
\mathrm{e}^{-
c\alpha(n (t+u))^{1/3}}\frac{\mathrm{d}u}{(t+u)^{1+\beta}} \biggr)^{1/\alpha} \biggr
\}^{\alpha} \\
&\le & Cn^{-(({\alpha+1})/{2} -
\beta)} I,
\end{eqnarray*}
where $\frac{\alpha+1}{2} - \beta>0$ and $I :=  \{\int_0^\infty
\mathrm{d}
t  (\int_0^\infty\mathrm{e}^{- c\alpha
(t+u)^{1/3}}(t+u)^{-1-\beta}\,\mathrm{d}u
 )^{1/\alpha}  \}^{\alpha} < \infty$.
Next,
\begin{eqnarray*}
I_{n2} &\le&Cn^{\alpha/2} \int_0^\infty
z^\beta\,\mathrm{d}z \biggl\{\int_{(0,4]^2}
\mathrm{e}^{- zt - c(nt)^{1/3} - c(\sqrt{n} |s|)^{1/2}}\,\mathrm{d}t \,\mathrm{d}s \biggr\}^{\alpha
}
\\
&\le& C \biggl\{\int_0^4
\mathrm{e}^{- c(nt)^{1/3}}\, \mathrm{d}t \biggl(\int_0^\infty
\mathrm{e}^{- \alpha
z t} z^\beta\,\mathrm{d}z \biggr)^{1/\alpha}
\biggr\}^{\alpha}
\\
&\le&C \biggl\{\int_0^\infty
\mathrm{e}^{- c(nt)^{1/3}} t^{-({1+\beta})/{\alpha
}} \,\mathrm{d}t \biggr\}^{\alpha}
\le Cn^{-(\alpha- 1 - \beta)} = o(1).
\end{eqnarray*}
Finally, using $\mathrm{e}^{- c(\sqrt{n} |s-v|)^{1/2}} \le\mathrm
{e}^{- (c/2)(\sqrt{n}
|v|)^{1/2}}$
for $|v| \ge2, |s|\le1 $ it easily follows $I_{n3} = O(\mathrm{e}^{-
c'n^{1/4}}) = o(1) (\exists c'>0)$, thus completing the proof of
\eqref{3NJ} and \eqref{fidi1} for $\gamma=\gamma_0 =1/2$.

\textit{Case} $\gamma> 1/2, 0< \beta< \alpha-1$. In this case,
\eqref{Jngamma} holds with
\begin{eqnarray}
&& G_{n\gamma}(u,v,z)
\nonumber
\\
&&\quad:=\sum_{j=1}^p \theta_j
\int_0^{{\lfloor nx_j\rfloor }/{n}} \,\mathrm{d}t n^{-1/2} \sum
_{s=1}^{\lfloor n^\gamma y_j\rfloor } g_3 \biggl(\lceil
nt\rceil -\lceil nu\rceil , s-\bigl\lceil n^\gamma v\bigr\rceil , 1 -
\frac{z}{n} \biggr) \mathbf{1} (0< z< n)
\nonumber
\\
\label{Gnplus1} &&\quad=\sum_{j=1}^p
\theta_j \int_0^{{\lfloor nx_j\rfloor }/{n}} \,\mathrm{d}t
\int_{\mathbb{R}} \,\mathrm{d}s \sqrt n g_3 \biggl(\lceil
nt\rceil -\lceil nu\rceil , \lceil \sqrt n s\rceil , 1 - \frac
{z}{n}
\biggr)
\\
&&\qquad {}\times\mathbf{1} \bigl(0< z< n, 1-\bigl\lceil n^\gamma v\bigr
\rceil \leq \lceil \sqrt n s\rceil \leq \bigl\lfloor n^\gamma
y_j\bigr\rfloor -\bigl\lceil n^\gamma v\bigr\rceil \bigr)
\nonumber
\\
&&\quad=:\sum_{j=1}^p \theta_j
\int_0^{x_j} \,\mathrm{d}t \int
_{\mathbb
{R}} \,\mathrm{d}s f_{nj}(t,s,u,v,z).
\nonumber
\end{eqnarray}
We first check the point-wise convergence: for any $ (u,z) \in\mathbb{R}
\times
\mathbb{R}_+,
v \in\mathbb{R}\setminus\{0, y_j\},  j=1,\ldots, p$
%
\begin{eqnarray}
G_{n\gamma}(u,v,z)  &\to &  G_\gamma(u,v,z)
\nonumber
\\[-8pt]
\label{konverg1}
\\[-8pt]
\nonumber
& :=& \sum
_{j=1}^p \theta_j \int
_0^{x_j} \,\mathrm{d}t \int_{\mathbb{R}}
\mathrm{d}s h_{3}(t-u, s, z) \mathbf{1}(0 < v < y_j),
\qquad n\to\infty.
\end{eqnarray}
To prove \eqref{konverg1}, note that from \eqref{green03}, \eqref
{limg3bdd} and $\gamma> 1/2$, for any $u < t \in\mathbb{R}$, $v \in
\mathbb{R}
\setminus\{0, y_j\},  j=1,\ldots, p$, $s\in\mathbb{R}$, and
$z>0$, we have
the point-wise convergences (as $n\to\infty$)
\begin{eqnarray*}
\sqrt n g_3 \biggl(\lceil nt\rceil -\lceil nu\rceil , \lceil \sqrt
n s\rceil , 1 - \frac
{z}{n} \biggr)\mathbf{1} (0< z< n) &\to&
h_{3}(t-u, s, z),
\\
\mathbf{1} \bigl(1-\bigl\lceil n^\gamma v\bigr\rceil \leq\lceil \sqrt n
s\rceil \leq \bigl\lfloor n^\gamma y_j\bigr\rfloor -\bigl
\lceil n^\gamma v\bigr\rceil \bigr) &\to& \mathbf{1}(0 < v <
y_j)
\end{eqnarray*}
and hence
%
\begin{equation}
\label{fn} f_{nj}(t,s;u,v,z) \to f_j(t,s;u,v,z) :=
h_{3}(t-u, s, z)\mathbf{1}(0 < v < y_j). 
\end{equation}
Using \eqref{fn}, relation \eqref{konverg1} can be shown similarly as
in the case $\gamma= \gamma_0$ above. Namely,
write $G_{n\gamma}(u,v,z)- G_{\gamma}(u,v,z)  =   \sum_{i=1}^3
\Gamma
_{ni}(u,v,z) $,
where, for $0< z < n $,
\begin{eqnarray*}
\Gamma_{n1}(u,v,z)&:=&\sum_{j=1}^p
\theta_j \int_0^{{\lfloor nx_j\rfloor }/{n}} \mathrm{d}t
\int_\mathbb{R}\bigl\{ f_{nj}(t,s;u,v,z)-
f_j(t,s;u,v,z)\bigr\} \mathbf{1}\bigl((t,s)\in D_j(
\epsilon_1)\bigr) \,\mathrm{d} s,
\\
\Gamma_{n2}(u,v,z)&:=&\sum_{j=1}^p
\theta_j \int_0^{{\lfloor nx_j\rfloor }/{n}} \mathrm{d}t
\int_\mathbb{R}f_j(t,s;u,v,z) \mathbf{1}
\bigl((t,s)\notin D_j(\epsilon _1)\bigr) \,\mathrm{d}s,
\\
\Gamma_{n3}(u,v,z)&:=&\sum_{j=1}^p
\theta_j \int_0^{\lfloor nx_j\rfloor/n} \mathrm{d}t
\int_\mathbb{R}f_{nj}(t,s;u,v,z) \mathbf{1}
\bigl((t,s)\notin D_j(\epsilon_1)\bigr) \,\mathrm{d}s,
\end{eqnarray*}
and where
$D_j(\epsilon_1) := \{ (t,s) \in(0,x_j]\times\mathbb{R}: t-u >
\epsilon_1,
|s-v| > \epsilon_1, |s| < 1/\epsilon_1 \} $.
Then \eqref{konverg1} follows if we show that, for any $ (u,z) \in
\mathbb{R}
\times\mathbb{R}_+,   v \in\mathbb{R}\setminus\{0, y\}$,
%
\begin{eqnarray}
 \lim_{n \to\infty}\bigl|\Gamma_{n1}(u,v,z)\bigr| &=& 0
\qquad \forall \epsilon_1 >0\quad \mbox{and}\nonumber
\\[-8pt]
\label{gamma21}
\\[-8pt]
\nonumber \lim
_{\epsilon_1 \to0} \limsup_{n\to\infty} \bigl|\Gamma
_{ni}(u,v,z) \bigr| &=& 0, \qquad i=2,3.
\end{eqnarray}
Here, the first relation in \eqref{gamma21} follows from the uniform
convergence statement of Lemma~\ref{lemma3N}, and
the second one from the dominating bound in \eqref{limg3bdd}; in particular,
\begin{eqnarray*}
&&\int_0^{x_j} \mathrm{d}t \int
_\mathbb {R}f_{nj}(t,s;u,v,z) \mathbf{1}(t-u \le
\epsilon_1) \,\mathrm{d}s
\\
&&\qquad\le\int_0^{\epsilon_1 + n^{-1}} \mathrm{d}t \int
_{\mathbb{R}} \biggl(\frac
{1}{\sqrt
{t}} \mathrm{e}^{- cs^2/t} +
\sqrt{n} \mathrm{e}^{- c(nt)^{1/3} -
c(\sqrt{n}|s|)^{1/2}} \biggr) \,\mathrm{d}s \le C\bigl(
\epsilon_1 + n^{-1}\bigr)
\end{eqnarray*}
vanishes as $n\to\infty$ and $\epsilon_1 \to0$.

With \eqref{konverg1} in mind, the convergence of integrals in \eqref
{gamma0} and \eqref{gamma2} can be established using
the dominated convergence theorem and the bound \eqref{limg3bdd} of
Lemma~\ref{lemma3N}, similarly
as in the case $\gamma= 1/2$ above.

 \textit{Case} $\gamma> 1/2, -(\alpha-1)/2 < \beta<0$.
 In this
case, \eqref{Jngamma} holds with $\chi_n(z) := (z/n^{2\gamma
})^{-\beta}
(\phi(1 - (z/n^{2\gamma}))/\phi_1)
\mathbf{1}(0< z < n^{2\gamma}) \to1 $ and
\begin{eqnarray*}
G_{n\gamma}(u,v,z) &:=&\sum_{j=1}^p
\theta_j \int_0^{{\lfloor nx_j\rfloor }/{n}} \mathrm{d}t
\int_0^{{\lfloor n^\gamma y_j\rfloor }/{n^\gamma}} \mathrm{d}s n^\gamma
\\
&&{}\times g_3 \biggl(\lceil nt\rceil -\bigl\lceil
n^{2\gamma} u\bigr\rceil , \bigl\lceil n^\gamma s\bigr\rceil -
\bigl\lceil n^\gamma v\bigr\rceil , 1 - \frac
{z}{n^{2\gamma}} \biggr)
\mathbf{1} \bigl(0< z< n^{2\gamma}\bigr)
\\
&=:&\sum_{j=1}^p \theta_j
\int_0^{{\lfloor nx_j\rfloor}/{n}} \mathrm{d}t \int
_0^{{\lfloor n^\gamma y_j\rfloor }/{n^\gamma}} \mathrm{d}s f_{n}(t,s;u,v,z).
\end{eqnarray*}
Note that in the above integral, variables $t$ and $u$ are rescaled by
$n $ and $n^{2\gamma} \gg n $, respectively.
Therefore, by \eqref{green03} the integrand
%
\begin{equation}
\label{fn1} f_{n}(t,s;u,v,z) \to f(s;u,v,z) :=
h_{3}(-u, s-v, z) \qquad \mbox{as } n\to\infty
\end{equation}
converges point-wise to $f(s;u,v,z)$ independent of $t$, for any $u<
0$, $s, v \in\mathbb{R}$, $s\in\mathbb{R}$, and $z>0$ fixed. By using
\eqref{fn1} and splitting $G_{n\gamma}(u,v,z)$ similarly as in the case
$\gamma= \gamma_0$ above, we can show
the uniform convergence in \eqref{gamma1} with
\begin{eqnarray*}
G_{\gamma}(u,v,z) &:=&\sum_{j=1}^p
\theta_j \int_0^{x_j} \mathrm{d}t
\int_0^{y_j} f(s;u,v,z)\, \mathrm{d} s = \sum
_{j=1}^p \theta_j
x_j \int_0^{y_j}
h_{3}(-u, s-v, z) \,\mathrm{d}s
\end{eqnarray*}
satisfying \eqref{Jgamma}; see the definition of
$F_{3\gamma}(x,y; u,v,z) = G_{\gamma}(u,v,z)$ in \eqref{Fgamma}. The
proof of \eqref{gamma2}
uses the dominating bound \eqref{limg3bdd} of Lemma~\ref{lemma3N}
similarly as in the previous cases.

\textit{Case} $0< \gamma< 1/2,  (\alpha-1)/2 < \beta< \alpha-1$.
We have \eqref{Jngamma} with
\begin{eqnarray*}
&&G_{n\gamma}(u,v,z)\\
&&\quad:=\sum_{j=1}^p
\theta_j \int_0^{{\lfloor nx_j\rfloor }/{n^{2\gamma}}} \mathrm{d}t
\int_0^{{\lfloor n^\gamma y_j\rfloor }/{n^\gamma}} \mathrm{d}s n^\gamma
g_3 \biggl(\bigl\lceil n^{2\gamma} t\bigr\rceil -\lceil nu
\rceil , \bigl\lceil n^\gamma s\bigr\rceil - \bigl\lceil
n^\gamma v\bigr\rceil , 1 - \frac{z}{n^{2\gamma}} \biggr)
\\
&&\quad=\sum_{j=1}^p \theta_j
\int_0^\infty\mathrm{d}w \int_0^{\lfloor n^\gamma y_j\rfloor /n^\gamma}
f_{nj}(w,s;u,v,z) \,\mathrm{d}s,
\end{eqnarray*}
where
\begin{eqnarray*}
f_{nj}(w,s;u,v,z)&:=&n^\gamma g_3 \biggl(\bigl
\lceil n^{2\gamma} w\bigr\rceil , \bigl\lceil n^\gamma s\bigr\rceil
- \bigl\lceil n^\gamma v\bigr\rceil , 1 - \frac{z}{n^{2\gamma}} \biggr)
\\
&& {}\times\mathbf{1} \biggl( \frac{1-\lceil n u\rceil }{n^{2\gamma}} < w < \frac{\lfloor
nx_j\rfloor - \lceil
nu\rceil
}{n^{2\gamma}}
\biggr)
\\
&\to&\mathbf{1}(0< u < x_j)h_3(w, s-v, z),\qquad n
\to\infty
\end{eqnarray*}
point-wise for each $u \in\mathbb{R}\setminus\{ 0, x_j\}, w > 0,
s\in
(0,y_j), v\in\mathbb{R}, s\ne v, z >0$ fixed, according
to Lemma~\ref{lemma3N}. This leads to the point-wise convergence of
integrals, namely,
%
\begin{eqnarray}
 G_{n\gamma}(u,v,z) &\to &  G_\gamma(u,v,z)
\nonumber
\\[-8pt]
\label{konverg2}
\\[-8pt]
\nonumber
&=& \sum
_{j=1}^p \theta_j
\mathbf{1}(0 < u < x_j) \int_0^\infty
\mathrm{d}w \int_0^{y_j} \mathrm{d}s
h_{3}(w,s-v, z),\qquad n\to\infty
\end{eqnarray}
similarly as in \eqref{konverg1} above. We omit the rest of the proof
of \eqref{gamma0} and \eqref{gamma2} which uses
\eqref{konverg2}, Lemma~\ref{lemma3N} and the dominated convergence theorem.

\textit{Case} $0< \gamma< 1/2,  -(\alpha-1)/2< \beta<
(\alpha -1)/2$.  We have \eqref{Jngamma} with
\[
G_{n\gamma}(u,v,z) := \sum_{j=1}^p
\theta_j \int_0^{{\lfloor nx_j\rfloor }/{n}} \mathrm{d}t
\int_0^{{\lfloor n^\gamma y_j\rfloor }/{n^\gamma}} f_{n}(t,s;u,v,z)
\,\mathrm{d}s
\]
and
\begin{eqnarray*}
f_{n}(t,s;u,v,z)&:=& n^{1/2} g_3 \biggl(\lceil
n t\rceil - \lceil nu\rceil , \bigl\lceil n^\gamma s\bigr\rceil - \bigl
\lceil n^{1/2} v\bigr\rceil , 1 - \frac{z}{n} \biggr) \mathbf{1}
(0< z \le n )
\\
&\to& h_3(t-u, -v, z)
\end{eqnarray*}
tending to a limit independent of $s $ for each $ t < u, s\in\mathbb
{R}, v\in
\mathbb{R}, z >0$ fixed, according
to Lemma~\ref{lemma3N} and using the fact that
$\sup_{s \in[0,y]} |\frac{\lceil n^\gamma s\rceil - \lceil
n^{1/2} v\rceil
}{n^{1/2}} - v | \to0$ for any $y>0$ as
$\gamma< 1/2 $. Whence, the point-wise convergence, as $n\to\infty$,
%
\begin{eqnarray}
G_{n\gamma}(u,v,z) &\to&  G_\gamma(u,v,z) = \sum
_{j=1}^p \theta_j
y_j \int_0^{x_j}
h_{3}(t-u,-v, z) \,\mathrm{d}t
\nonumber
\\[-8pt]
\label{konverg3}
\\[-8pt]
\nonumber
&=& \sum_{j=1}^p
\theta_j F_{3\gamma
}(x_j,y_j;u,v,z),
\end{eqnarray}
can be obtained.
The details of the proof of \eqref{konverg3} and subsequently \eqref
{gamma1} and \eqref{gamma2} are similar as in other cases
above.

This proves \eqref{fidi1}, and hence the limit in \eqref{3Nlim0} in all
cases of $\gamma$ and $\beta$ under consideration.
The second statement of the theorem follows from \eqref{3Nlim0} and
Proposition~\ref{3Nexist}.
Theorem~\ref{3Ntheo} is proved.
\end{pf*}

\section{Scaling transition in the aggregated 4N model}\label{sec4}

In this section, we discuss scaling transition and Type I
isotropic distributional LRD property for the aggregated 4N model
${\mathfrak X}_4$ in \eqref{aggre34}.
Recall that $g_4(t,s,a)$ in \eqref{aggre34} is the Green function of
the random walk
$\{W_k \} $ on $\mathbb{Z}^2 $ with one-step transition probabilities
shown in
Figure~\ref{fig1}(b). Recall that
\[
h_4(t,s,z) = \frac{2}{\pi}K_0 \bigl(2\sqrt{z
\bigl(t^2 + s^2\bigr)} \bigr) = \frac
{2}{\pi} \int
_0^\infty w^{-1} \mathrm{e}^{- zw - (t^2+s^2)/w }
\,\mathrm{d}w, \qquad (t,s) \in\mathbb{R} ^2_0, z >0
\]
is the potential of the Brownian motion in $\mathbb{R}^2$ with covariance
matrix $\operatorname{diag}(1/2, 1/2)$,
written
via $K_0$, the modified Bessel function of second kind. See \cite
{ito1965}, Chapter~7.2.

For any $\gamma>0$,
introduce a RF $V_{4\gamma} = \{V_{4 \gamma}(x,y); (x,y) \in\bar
\mathbb{R}
^2_+\}$ as a stochastic integral
%
\begin{eqnarray}
\label{4Lfield0} V_{4 \gamma}(x,y)&:=&\int_{\mathbb{R}^2 \times\mathbb{R}_+}
F_{4\gamma}(x,y; u,v,z) \mathcal{M}(\mathrm{d}u, \mathrm{d}v, \mathrm{d}z),
\end{eqnarray}
where $F_{4\gamma}(x,y; u,v,z)$ is defined as
%
\begin{eqnarray}
\label{F4gamma} F_{4\gamma}& :=& %
\cases{ \displaystyle \int
_0^x \!\int_0^y
h_4(t-u, s-v, z) \,\mathrm{d}t \,\mathrm{d}s, &\quad $\gamma = 1$,
\vspace*{3pt}
\cr
\displaystyle \mathbf{1}(0<v<y) \int_0^x
\mathrm{d}t \int_\mathbb{R}h_4(t-u, w, z)
\,\mathrm{d}w, & \quad$\gamma> 1, \beta> (\alpha-1)/2$,\vspace*{3pt}
\cr
\displaystyle \mathbf{1}(0<u<x) \int_\mathbb{R}\mathrm{d}w \int
_0^y h_4(w, s-v, z)\, \mathrm{d}s, &
\quad $\gamma< 1, \beta> (\alpha-1)/2$, \vspace*{3pt}
\cr
\displaystyle x \int
_0^y h_4(u, s-v, z) \,\mathrm{d}s, &
\quad $\gamma>1, 0< \beta< (\alpha-1)/2$,\vspace*{3pt}
\cr
y \int
_0^x h_4(t-u, v, z) \,\mathrm{d}t, &
\quad $\gamma< 1, 0< \beta< (\alpha-1)/2$,}
\end{eqnarray}
and where
$\mathcal{M}$ is the same $\alpha$-stable random measure on $\mathbb{R}^2
\times\mathbb{R}_+$ as in
\eqref{3Lfield0}.

\begin{proposition} \label{4Nexist}
\textup{(i)} $V_{4\gamma}$
in (\ref{4Lfield0}) is well-defined for any $\gamma>0, 1< \alpha\le
2, 0< \beta< \alpha-1$
with exception of $\gamma\ne1, \beta= (\alpha-1)/2$.
It has $\alpha$-stable finite-dimensional distributions and stationary
rectangular increments
in the sense of \eqref{Vstatinc}.

\textup{(ii)} $V_{4\gamma}$ is OSRF: for any $\lambda>0$,
$\{V_{4\gamma}(\lambda x, \lambda^\gamma y); (x,y) \in\mathbb
{R}^2_+\}
\eqfdd
\{\lambda^{H(\gamma)} V_{4\gamma}(x,y); (x,y) \in\mathbb{R}^2_+\}$,
with
%
\begin{equation}
\label{Hgamma4N} H(\gamma):= %
\cases{ \displaystyle \frac{2(\alpha- \beta)}{\alpha}, &
\quad $\gamma= 1$, \vspace*{3pt}
\cr
\displaystyle\frac{\gamma-1 + 2(\alpha- \beta)}{\alpha}, &\quad $
\gamma> 1, \beta> (\alpha-1)/2$,\vspace*{3pt}
\cr
\displaystyle\frac{\alpha+ \alpha\gamma- 2\beta\gamma}{\alpha},
& \quad$\gamma> 1, \beta< (\alpha-1)/2$, \vspace*{3pt}
\cr
\displaystyle
\frac{1-\gamma+ 2\gamma(\alpha- \beta)}{\alpha}, & \quad$\gamma< 1, \beta> (\alpha-1)/2$,\vspace*{3pt}
\cr
\displaystyle \frac{\alpha+ \alpha\gamma- 2\beta}{\alpha}, &\quad$\gamma< 1, \beta < (\alpha-1)/2$.}
\end{equation}

\textup{(iii)} RFs $V_{4\gamma} = V_{4,+}  (\gamma>1) $ and
$V_{4\gamma}
= V_{4,-}  (\gamma<1) $
do not depend on $\gamma$ for $\gamma> 1$ and $\gamma< 1$.

\textup{(iv)} RF $V_{4\gamma}$ has properly dependent rectangular
increments for $\gamma= 1$ and does not have properly
dependent rectangular increments for $\gamma\ne1 $.

\textup{(v)} For $\alpha= 2$, the RFs
%
\begin{eqnarray}
V_{4,+} &\eqfdd& \kappa_{4,+} %
\cases{
B_{(3/2)- \beta,1/2}, & \quad $1/2< \beta< 1$, \vspace*{3pt}
\cr
B_{1,1 - \beta}, &
\quad $0 < \beta< 1/2$,}
\nonumber
\\[-8pt]
\\[-8pt]
\nonumber
V_{4,-} &\eqfdd& \kappa_{4,-} %
\cases{B_{1/2, (3/2)- \beta}, &\quad $1/2< \beta< 1$, \vspace*{3pt}
\cr
B_{1-\beta,1}, & \quad $0 < \beta< 1/2$, }
\end{eqnarray}
agree, up to some constants $\kappa_{4,\pm} = \kappa_{4,\pm} (\beta)
\neq0$, with fractional Brownian sheet $B_{H_1,H_2}$
where one of the parameters $H_1, H_2$ equals $1/2$ or 1.
\end{proposition}

\begin{pf}
(i) As in the proof of Proposition~\ref{3Nexist}(i), we
show $J_\gamma:= \int_{\mathbb{R}^2 \times\mathbb{R}_+}
(F_{4\gamma}(1,1;
u,v, z))^\alpha\,\mathrm{d}\mu< \infty$ only.
First, consider the case $\gamma= 1$. We have
$J_1 = C\int_{\mathbb{R}^2 \times\mathbb{R}_+}
 ( \int_{(0,1]^2} K_0(2\sqrt{z} \|v-w\|) \,\mathrm{d}v
)^\alpha
z^\beta\,\mathrm{d}w \,\mathrm{d}z
< \infty$. Here, $\|x\|^2:=x_1^2+x_2^2$, for $x=(x_1, x_2)\in\mathbb{R}^2$.
Split $J_1 = J' + J''$, where
$J' := \int_{ \{\|w\| \le\sqrt{2} \} \times\mathbb{R}_+ } \cdots,
 J'' :=
\int_{ \{\|w\| > \sqrt{2} \} \times\mathbb{R}_+ } \cdots.$
By Minkowski's inequality,\vspace*{-2pt}
\begin{eqnarray*}
J''&\le& C \biggl\{ \int_{\{\|v\| \le\sqrt{2} \}}
\mathrm{d}v \biggl[ \int_{
\{\|
w\| >\sqrt{2} \} \times\mathbb{R}_+ } K^\alpha_0\bigl(2
\sqrt{z} \|v-w\|\bigr) z^\beta\,\mathrm{d}z \,\mathrm{d}w \biggr]^{1/\alpha}
\biggr\}^\alpha
\\[-2pt]
&\le&C \biggl\{ \int_{\{\|v\| \le\sqrt{2} \}} \,\mathrm{d}v \biggl[ \int
_{ \{
\|w\|
> \sqrt{2} \}} \|v-w\|^{-2 - 2\beta}\, \mathrm{d}w
\biggr]^{1/\alpha} \biggr\}^\alpha
\\[-2pt]
&\le&C \biggl\{ \int_{\{\|v\| \le\sqrt{2} \}}\bigl(\sqrt{2}- \|v\| \bigr)^{-2\beta
/\alpha}
\,\mathrm{d}v \biggr\}^\alpha < \infty,
\end{eqnarray*}
where we used the facts that $\int_{0}^\infty
K^\alpha_0(2\sqrt{z}) z^\beta\,\mathrm{d}z < \infty$ and $0 < \beta
< \alpha
-1 \le2 $. Next,
\begin{eqnarray*}
J'&\le&C \int_{\{\|w\| \le\sqrt{2} \}} \mathrm{d}w \int
_0^\infty z^\beta \,\mathrm{d}z \biggl(
\int_{ \{\|v\| \le\sqrt{2} \} } K _0\bigl(2\sqrt{z} \|v\|\bigr) \,\mathrm{d}v
\biggr)^\alpha
\\[-2pt]
&\le&C\int_0^\infty z^\beta\,\mathrm{d}z
\biggl( \int_{0}^{\sqrt{2}} K _0(2
\sqrt{z} r) r\, \mathrm{d}r \biggr)^\alpha
\\[-2pt]
&\le& C\int_0^\infty z^\beta
\bigl(z^{-\alpha/2} \mathbf{1} (0<z < 1) + z^{-\alpha} \mathbf{1}(z\ge1)
\bigr) \,\mathrm{d}z < \infty,
\end{eqnarray*}
where we used $0< \beta< \alpha- 1$ and the\vspace*{-3pt} inequality
\[
\int_{0}^{\sqrt{2}} K _0(2\sqrt{z} r) r
\,\mathrm{d}r \le C %
\cases{ z^{-1/2}, &\quad $0<z \le 1$,
\vspace*{2pt}
\cr
z^{-1}, &\quad $z> 1$,}
\]
which is a consequence of the fact that the function $r \mapsto r
K_0(r) $ is bounded and integrable on $(0,\infty)$.
This proves $J_1 < \infty$.

Next, let $\gamma> 1, (\alpha-1)/2 < \beta< \alpha-1$. Using
$h_{4\star} (u,z) := \int_\mathbb{R}h_4(u, w, z) \,\mathrm{d}w
= \frac{2}{\pi} \int_\mathbb{R}K_0 (2\times\sqrt{z(u^2 + w^2)} )
\,\mathrm{d}w =
\frac
{2}{\pi} \sqrt{\frac{u}{4 z^{1/2}}} K_{-1/2}(2\sqrt{z} |u|) =
\sqrt{\frac{1}{4 \pi z}} \mathrm{e}^{- 2 \sqrt{z} |u|} $
(\cite{gra1962}, 6.596, 8.469), we obtain
$J_\gamma\le C\int_{\mathbb{R}} \mathrm{d}u \int_{\mathbb{R}_+}
z^\beta\,\mathrm{d}z  ( \int_0^1
h_{4 \star}(t-u,z)\, \mathrm{d}t  )^\alpha
\le C  \{\int_{|u| \le2} \cdots+ \int_{|u| > 2} \cdots \}
=: C
\{J'_\gamma+ J''_\gamma\}$,\vspace*{-2pt}
where
\begin{eqnarray*}
J'_\gamma &\le& C \int_0^\infty
z^\beta\,\mathrm{d}z \biggl( \int_0^1
h_{4
\star
}(t,z) \,\mathrm{d}t \biggr)^\alpha \le C \int
_0^\infty z^{\beta- (\alpha/2)}\, \mathrm{d}z \biggl(
\int_0^1 \mathrm{e}^{-
2\sqrt{z} t}\, \mathrm{d}t
\biggr)^\alpha
\\[-2pt]
&\le& C \int_0^\infty z^{\beta- \alpha}\,
\mathrm{d}z \bigl(1- \mathrm {e}^{- 2\sqrt
{z}} \bigr)^\alpha,
\end{eqnarray*}
where the last integral converges for any $0< \beta< \alpha-1, 1<
\alpha\le2$. Next,\vspace*{-2pt}
\[
J''_\gamma\le C\int_1^\infty
\mathrm{d}u \int_0^\infty z^{\beta-
(\alpha/2)}
\mathrm{e}^{- 2\sqrt{z}u}\, \mathrm{d}z \le C\int_0^\infty
z^{\beta- (1+ \alpha)/2} \mathrm{e}^{- 2 z} \,\mathrm{d}z < \infty
\]
provided $\beta> (\alpha-1)/2 $ holds. Hence, $J_\gamma< \infty$.

Consider $J_\gamma$ for $\gamma> 1, 0< \beta< (\alpha-1)/2$. We have
$J_\gamma\le C\int_{\mathbb{R}} \mathrm{d}u \int_\mathbb
{R}\,\mathrm{d}v \int_{\mathbb{R}_+} z^\beta\,\mathrm{d}
z
( \int_0^1 h_{4}(u, s-v,z) \,\mathrm{d}s  )^\alpha
\le C  \{\int_{|v| \le2} \cdots+ \int_{|v| > 2} \cdots \}
=: C
\{J'_\gamma+ J''_\gamma\}$.
By Minkowski's inequality,
\begin{eqnarray*}
J'_\gamma &\le&C\int_0^\infty
\mathrm{d}u \int_0^\infty z^\beta
\mathrm{d}z \biggl( \int_0^1
h_{4}(u,s,z) \,\mathrm{d}s \biggr)^\alpha
\\
&\le& C \biggl\{ \int_0^1 \mathrm{d}s \biggl[
\int_0^\infty\mathrm {d}u \int
_0^\infty z^\beta K^\alpha_0
\bigl(2 \sqrt{z\bigl(t^2 + u^2\bigr)}\bigr) \,\mathrm{d}z
\biggr]^{1/\alpha} \biggr\} ^\alpha
\\
&\le&C \biggl\{ \int_0^1 \mathrm{d}s \biggl[
\int_0^\infty\frac
{\mathrm{d}u}{(t^2 +
u^2)^{\beta+1}}
\biggr]^{1/\alpha} \biggr\}^\alpha \le C \biggl\{ \int
_0^1 \mathrm{d}s \biggl[\frac{1}{s^{2\beta+1}}
\biggr]^{1/\alpha
} \biggr\}^\alpha < \infty
\end{eqnarray*}
since $\beta< (\alpha-1)/2$. Next,
\begin{eqnarray*}
J''_\gamma &\le&C \int_1^\infty
\mathrm{d}v \int_0^\infty\mathrm{d}u \int
_0^\infty z^\beta\, \mathrm{d}z
h^\alpha_{4}(u,v,z)
\\
&\le& C\int_1^\infty\mathrm{d}v \int
_0^\infty\frac{\mathrm
{d}u}{(u^2 +
v^2)^{\beta
+1}} \le C\int
_1^\infty\frac{\mathrm{d}v}{v^{1+2\beta}} < \infty.
\end{eqnarray*}
Hence, $J_\gamma< \infty$ for $\gamma> 1$. The case $0< \gamma< 1 $
follows by symmetry. This proves
the existence of $V_{4\gamma} $
for all choices of $\alpha, \beta, \gamma$ in \eqref{F4gamma}. The
remaining facts in (i) are similar
as in Proposition~\ref{3Nexist}.

(ii) Follows analogously as in Proposition~\ref{3Nexist}(ii).

(iii) Follows from the definition of the integrand
$F_{4\gamma}$
in \eqref{F4gamma}.

(iv) The proof is completely similar to that of Proposition~\ref{3Nexist}(iii),
taking into account the form of $ V_{4\gamma}$
in \eqref{4Lfield0} and the fact that $h_4(u,v,z) $ is everywhere
positive on $\mathbb{R}^2 \times\mathbb{R}_+$.

(v) Follows from the OSRF property in (ii) analogously as in
Proposition~\ref{3Nexist}(v).
Proposition~\ref{4Nexist} is proved.
\end{pf}

The main result of this section is Theorem~\ref{4Ntheo}. Its proof is
based on the asymptotics of the Green
function $g_4$ in Lemma~\ref{lemma4N}, below. The proof of Lemma~\ref{lemma4N}
can be found at \surl{http://arxiv.org/abs/1303.2209v3}.

\begin{lem} \label{lemma4N}
For any $ (t,s,z ) \in\mathbb{R}^2_0 \times(0,\infty)$
%
\begin{equation}
\label{limg4} \lim_{\lambda\to\infty} g_4 \biggl([\lambda
t], [\lambda s], 1- \frac
{z}{\lambda^2} \biggr) = h_4(t,s,z) =
\frac{2}{\pi} K_0 \bigl(2 \sqrt{z\bigl(t^2 +
s^2\bigr)} \bigr).
\end{equation}
The convergence in \eqref{limg4} is uniform on any relatively compact
set $\{\epsilon< |t|+|s| < 1/\epsilon\}
\times\{\epsilon< z < 1/\epsilon\}
\subset\mathbb{R}^2_0 \times\mathbb{R}_+,  \epsilon>0$.

Moreover, there exists constants $C, c >0$ such that for
all sufficiently large $\lambda$ and
any $(t,s,z) \in\mathbb{R}^2_0 \times(0, \lambda^2)$ the following
inequality holds:
%
\begin{equation}
\label{limg4bdd} g_4 \biggl([\lambda t], [\lambda s], 1-
\frac{z}{\lambda^2} \biggr) < C \bigl\{h_4(t,s, z) +
\mathrm{e}^{- c \sqrt{\lambda} (|t|^{1/2} + |s|^{1/2})} \bigr\}.
\end{equation}
\end{lem}

\begin{theorem} \label{4Ntheo}
Assume that the mixing density $\phi$ is
bounded on $[0,1)$ and satisfies (\ref{mixdensity}), where
%
\begin{equation}
\label{4bcond} 0 < \beta< \alpha- 1, \qquad 1< \alpha\le2, \qquad \beta\neq (
\alpha-1)/2.
\end{equation}
Let $\mathfrak{X}_4$ be the aggregated 4N model in \eqref{aggre34}.
Then for any $\gamma>0$
%
\begin{eqnarray}
\label{4Nlim} && n^{-H(\gamma)} \sum_{t=1}^{[nx]}
\sum_{s=1}^{[n^\gamma y]} \mathfrak
{X}_4(t,s) \limfdd V_{4\gamma}(x,y), \qquad x, y >0, n \to
\infty,
\end{eqnarray}
where $H(\gamma)$ and $V_{4\gamma}$ are given in \eqref{Hgamma4N} and
\eqref{4Lfield0}, respectively.
As a consequence,
the RF  $\mathfrak{X}_4$ exhibits scaling transition at $\gamma_0 =
1$ and
enjoys
Type I isotropic distributional LRD property
in the sense of Definition~\ref{Adlm}.
\end{theorem}

\begin{pf}
Similarly, as in the proof of Theorem~\ref{3Ntheo}, it suffices to prove
the limit
%
\begin{equation}
\label{lim4N} \lim_{n \to\infty} J_{n\gamma} =
J_{\gamma},
\end{equation}
where
%
\begin{eqnarray}
 J_{n \gamma} &:=&n^{-\alpha H(\gamma)}\sum
_{(u,v) \in\mathbb{Z}^2} \mathrm {E} \Biggl| \sum_{j=1}^p
\theta_j \sum_{1\le t \le[nx_j], 1\le s \le[n^\gamma y_j]}
g_4(t-u, s-v, A) \Biggr|^\alpha,
\nonumber
\\[-8pt]
\label{Jconv4}
\\[-8pt]
\nonumber
J_\gamma&:=&\int_{\mathbb{R}^2 \times\mathbb{R}_+} \bigl|G_\gamma
(u,v,z)\bigr|^\alpha\,\mathrm{d} \mu,\qquad G_\gamma(u,v,z) := \sum
_{j=1}^p \theta_j
F_{4\gamma}(x_j,y_j; u,v,z),
\nonumber
\end{eqnarray}
for any $p\in\mathbb{N}_+, \theta_j \in\mathbb{R}, (x_j,y_j) \in
\mathbb{R}^2_+, j=1,
\ldots, p$.
The proof of \eqref{lim4N} follows the same strategy as in
the case of Theorem~\ref{3Ntheo}, that is, we write $J_{n\gamma} $ as a
Riemann sum
approximation
%
\begin{equation}
\label{Jngamma4} J_{n\gamma} = \int_{\mathbb{R}^2 \times\mathbb{R}_+}
\bigl|G_{n\gamma
}(u,v,z)\bigr|^\alpha \chi _n(z) \mu(
\mathrm{d}u, \mathrm{d}v, \mathrm{d}z),
\end{equation}
to the integral $J_\gamma$, where $\chi_n(z) \to1  (n \to\infty)$
boundedly in $z >0$, and
$G_{n\gamma}: \mathbb{R}^2 \times\mathbb{R}_+ \to\mathbb{R}$ are
some functions
tending to $G_\gamma$ in \eqref{Jconv4}. We use Lemma~\ref{lemma4N}
and the dominated convergence
theorem to deduce the convergence in \eqref{lim4N}.
Because of the differences in the form of the integrand in \eqref
{F4gamma}, several cases of $\gamma$ and $\beta$
need to be discussed separately.
The approximation is similar as in the proof of Theorem~\ref{3Ntheo}
and is discussed briefly below.

For $\epsilon>0$, denote $W_\epsilon:= \{ (u,v,z) \in\mathbb{R}^2
\times\mathbb{R}
_+: |u|+|v|< 1/\epsilon,  \epsilon< z < 1/\epsilon\},
W_\epsilon^c := (\mathbb{R}^2 \times\mathbb{R}_+)\setminus
W_\epsilon$. Similarly, as
in Theorem~\ref{3Ntheo}, \eqref{lim4N} follows from
%
\begin{equation}
\label{gamma40} \lim_{n \to\infty} \int_{W_\epsilon}
\bigl|G_{n\gamma}(u,v,z) - G_\gamma (u,v,z)\bigr|^\alpha\,\mathrm{d}
\mu = 0\qquad \forall \epsilon>0,
\end{equation}
and
%
\begin{equation}
\label{gamma42} \lim_{\epsilon\to0} \limsup_{n \to\infty}
\int_{W_\epsilon^c} \bigl|G_{n\gamma}(u,v,z)\bigr|^\alpha
\,\mathrm{d}\mu = 0.
\end{equation}

\textit{Case} $\gamma= \gamma_0 = 1$.
In this case, \eqref{Jconv4}
and \eqref{Jngamma4} hold with
$G_{\gamma_0}(u,v,z):= \sum_{j=1}^p \theta_j \times\int_0^{x_j} \!\int_0^{y_j}
h_4(t-u, s-v, z) \,\mathrm{d}t\, \mathrm{d}s$ and
\begin{eqnarray*}
G_{n\gamma_0}(u,v,z)&:=&\sum_{j=1}^p
\theta_j \int_0^{\lfloor nx_j\rfloor/n}\! \int
_0^{\lfloor ny_j\rfloor/n} g_4 \biggl(\lceil nt\rceil
-\lceil nu\rceil , \lceil ns\rceil -\lceil nv\rceil , 1 - \frac
{z}{n^2}
\biggr)
\\
&&{}\times\mathbf{1}\bigl(0<z<n^2\bigr) \,\mathrm{d}t\, \mathrm{d}s.
\end{eqnarray*}
Then, by splitting $G_{n\gamma_0}(u,v,z) - G_{\gamma_0}(u,v,z) = \sum_{i=1}^3 \Gamma_{ni}(u,v,z)$ and using Lemma~\ref{lemma4N}
similarly as in the proof of Theorem~\ref{3Ntheo}, Case $\gamma= 1/2$,
relation \eqref{gamma40} can be obtained.

Consider \eqref{gamma42}. Since $G_{\gamma_0} \in L^\alpha(\mu)$, see
the proof of Proposition~\ref{4Nexist}(i),
relation \eqref{gamma42} holds with $G_{n\gamma_0}$ replaced by
$G_{\gamma_0}$. Hence and with \eqref{limg4bdd} in mind, it
suffices to check \eqref{gamma42} with $G_{n\gamma_0}$ replaced by
$\widetilde G_n(u,v,z) := \mathbf{1}(0< z < n^2) \int_0^{1}\! \int_0^{1}
\mathrm{e}^{-c
(\sqrt{n|t-u|} + \sqrt{n|s-v|})} \,\mathrm{d}t\, \mathrm{d}s$, which
follows from
%
\begin{eqnarray}
\label{4NJ} \widetilde J_n&:=&\int_{\mathbb{R}^2 \times\mathbb{R}_+}
\bigl(\widetilde G_n(u,v,z)\bigr)^\alpha\,\mathrm{d}\mu = O
\bigl(n^{2(\beta- \alpha+1)}\bigr) = o(1).
\end{eqnarray}
We have
$\widetilde J_n \le C n^{2\beta+ 2}  \{ \int_{\mathbb{R}}
(\int_{0}^1 \mathrm{e}
^{-c \sqrt{n|t-u|}}\,\mathrm{d}t  )^\alpha\,\mathrm{d}u  \}^2$, where\vspace*{1pt}
$\int_{\mathbb{R}}  (\int_{0}^1 \mathrm{e}^{-c \sqrt
{n|t-u|}}\,\mathrm{d}t
)^\alpha\,\mathrm{d}u
\le \int_{\{|u|< 2 \}}
(\cdots)^\alpha\,\mathrm{d}u + \int_{\{|u|\ge2 \}}
(\cdots)^\alpha\,\mathrm{d}u  =: i'_n + i''_n$.  Here,\vspace*{1.5pt} $i'_n \le C
(\int_0^3 \mathrm{e}^{-c \sqrt{nv}}\, \mathrm{d}v  )^\alpha
\le C/n^\alpha$ and $i''_n \le C\int_2^\infty\mathrm{e}^{- c \alpha
\sqrt
{n(u-1)}} \,\mathrm{d}u = O(\mathrm{e}^{- c' \sqrt{n}}),  c' >0$.
This proves \eqref{4NJ} and \eqref{gamma42}.

\textit{Case} $\gamma> 1, (\alpha-1)/2 < \beta< \alpha-1$.
We have
\eqref{Jngamma4} with $G_\gamma(u,v,z)
= \sum_{j=1}^p \theta_j \mathbf{1}(0< v <y_j) \int_0^{x_j} \mathrm{d}t
\int_{\mathbb{R}}
h_4(t-u,s,z)\, \mathrm{d}s $ and
\begin{eqnarray*}
G_{n\gamma}(u,v,z)&:=&\sum_{j=1}^p
\theta_j \int_0^{\lfloor nx_j\rfloor/n} \mathrm{d}t
n^{-1} \sum_{s=1}^{\lfloor n^\gamma y_j\rfloor }
g_4 \biggl(\lceil nt\rceil - \lceil nu\rceil , s-\bigl\lceil
n^\gamma v\bigr\rceil , 1 - \frac{z}{n^2} \biggr)
\\
&& {}\times\mathbf{1} \bigl(0< z< n^2\bigr)
\\
&=&\sum_{j=1}^p \theta_j
\int_0^{\lfloor nx_j\rfloor/n} \mathrm{d}t \int
_{\mathbb{R}} g_4 \biggl(\lceil nt\rceil -\lceil nu
\rceil , \lceil ns\rceil , 1 - \frac
{z}{n^2} \biggr)
\\
&& {}\times\mathbf{1} \bigl(0< z< n^2, 1-\bigl\lceil
n^\gamma v\bigr\rceil \leq\lceil ns\rceil \leq\bigl\lfloor
n^\gamma y_j\bigr\rfloor -\bigl\lceil n^\gamma v
\bigr\rceil \bigr)\, \mathrm {d}s
\\
&=:&\sum_{j=1}^p \theta_j
\int_0^{x_j} \mathrm{d}t \int
_{\mathbb
{R}} f_{nj}(t,s;u,v,z) \,\mathrm{d}s,
\end{eqnarray*}
cf. \eqref{Gnplus1}.
From \eqref{green04}, \eqref{limg4bdd} and $\gamma> 1$, for any $u, t
\in\mathbb{R}$, $u\neq t$, $v \in\mathbb{R}\setminus\{0,y_j\},
j=1, \ldots, p$,
$s$, and $z>0$, we have point-wise convergences
\begin{eqnarray*}
g_4 \biggl(\lceil nt\rceil -\lceil nu\rceil , \lceil n\rceil , 1 -
\frac{z}{n^2} \biggr)\mathbf{1} \bigl(0< z< n^2\bigr) &\to&
h_{4}(t-u, s, z),
\\
\mathbf{1} \bigl(1-\bigl\lceil n^\gamma v\bigr\rceil \leq\lceil ns
\rceil \leq \bigl\lfloor n^\gamma y_j\bigr\rfloor -\bigl
\lceil n^\gamma v\bigr\rceil \bigr) &\to& \mathbf{1}(0 < v <
y_j)
\end{eqnarray*}
implying $f_{nj}(t,s,u,v,z) \to  h_{4}(t-u, s, z)\mathbf{1}(0 < v < y_j)$
similarly as in \eqref{fn} in
the proof of
Theorem~\ref{3Ntheo}. The remaining details of the proof of \eqref
{gamma40} and
\eqref{gamma42} are similar as in Theorem~\ref{3Ntheo}, Case $\gamma>
1/2, 0 < \beta< \alpha-1$.

\textit{Case} $\gamma> 1, 0< \beta< (\alpha-1)/2$.  We have
\eqref{Jngamma4} with
$G_\gamma(u,v,z)
= \sum_{j=1}^p \theta_j x_j \int_0^{y_j} h_4(-u,s-v,z) \,\mathrm{d}s$ and
\begin{eqnarray*}
G_{n\gamma}(u,v,z) &:=&\sum_{j=1}^p
\theta_j \int_0^{\lfloor nx_j\rfloor/n} \,\mathrm{d}t
\int_0^{\lfloor n^\gamma y_j\rfloor/n^\gamma} g_4 \biggl(\lceil nt
\rceil - \bigl\lceil n^{\gamma} u\bigr\rceil , \bigl\lceil n^\gamma
s\bigr\rceil - \bigl\lceil n^\gamma v\bigr\rceil , 1 -
\frac{z}{n^{2\gamma}} \biggr)
\\
&&{}\times\mathbf{1} \bigl(0< z< n^{2\gamma}\bigr)
\\
&=:&\sum_{j=1}^p \theta_j
\int_0^{\lfloor nx_j\rfloor/n} \mathrm{d}t \int
_0^{\lfloor n^\gamma y_j\rfloor/n^\gamma} f_n(t,s;u,v,z)\,
\mathrm{d}s,
\end{eqnarray*}
where $f_n(t,s;u,v,z) \to f(s;u,v,z) := h_{4}(-u, s-v, z) $ tends to a
limit independent of $t$, as $n\to\infty$. Again, we omit the details
of the proof of \eqref{gamma40} and
\eqref{gamma42} which are similar as in Theorem~\ref{3Ntheo}, Case
$\gamma> 1/2, -(\alpha-1)/2 < \beta< (\alpha-1)/2$.

Case $0<\gamma< 1$ in \eqref{4Nlim}
follows from case $\gamma> 1$ by lattice isotropy of the 4N model.
This ends the proof of \eqref{4Nlim}. The second statement
of the theorem follows from Proposition~\ref{4Nexist}.
Theorem~\ref{4Ntheo} is proved. \hfill$\Box$

The following proposition obtains the asymptotic behavior of the
covariance function $r_4(t,s) = \mathrm{E}\mathfrak{X}_4 (t,s)
\mathfrak
{X}_4(0,0) $
of the aggregated Gaussian RF $\mathfrak{X}_4$ in \eqref{aggre34}
($\alpha= 2$).
The proof of Proposition~\ref{4Ncov} uses Lemma~\ref{lemma4N} and is omitted.
\end{pf}

\begin{proposition} \label{4Ncov}
Assume $\alpha=2$ and the conditions
of Theorem~\ref{4Ntheo}.
Then for any $(t,s) \in\mathbb{R}^2_0$
%
\begin{eqnarray}
\label{4Ncovasy} \lim_{\lambda\to\infty} \lambda^{2\beta}
r_4\bigl([\lambda t], [\lambda s]\bigr)&=& \frac{\sigma^2 \phi_1\Gamma(\beta+1)\Gamma(\beta)}{\pi}
\bigl(t^2+s^2\bigr)^{-\beta}.
\end{eqnarray}
\end{proposition}

\section{Auxiliary results}\label{sec5}

This section obtains conditions for the existence of a stationary
solution of a general
random-coefficient nearest-neighbor autoregressive RF in \eqref{NN}.
We also
discuss contemporaneous aggregation
of \eqref{NN} under the assumption that the innovations
belong the domain of attraction of $\alpha$-stable law,
$0< \alpha\le2$.

\subsection{Existence of random-coefficient autoregressive RF}

Consider a general random-coefficient nearest-neighbor autoregressive
RF on $\mathbb{Z}^2$:
%
\begin{eqnarray}
\label{NN} X(t,s)&=&\sum_{|u|+|v|=1} a(u,v) X(t+u,s+v)
+ \varepsilon (t,s), \qquad (t,s) \in\mathbb{Z}^2,
\end{eqnarray}
where $\{\varepsilon (t,s); (t,s) \in\mathbb{Z}^2 \} $ are i.i.d.
r.v.'s with finite
$p$th moment, $p \in(0,2]$,
and $a(t,s) \ge0, |t|+|s|=1 $ are \textit{random} coefficients
independent of $\{\varepsilon (t,s)\}$ and satisfying
%
\begin{equation}
\label{Acond} A := \sum_{|t|+|s| =1} a(t,s) \in (0,1)
\qquad \mbox {a.s.}
\end{equation}
Set also $a(t,s) := 0, (t,s) \in\mathbb{Z}^2, |t|+|s|\neq1$.
Clearly, the 3N and 4N models in \eqref{3N} and \eqref{4N} are
particular cases of \eqref{NN}.

Let us discuss solvability of \eqref{NN}. We will show that under
certain conditions
this equation
admits a stationary solution given by the convergent series
%
\begin{eqnarray}
\label{stationary0} X(t,s)&=&\sum_{(u,v) \in\mathbb{Z}^2 } g(t-u,s-v,
\mathbf{a})\varepsilon (u,v), \qquad (t,s) \in\mathbb{Z}^2,
\end{eqnarray}
where $g(t,s,\mathbf{a}),  (t,s) \in\mathbb{Z}^2, \mathbf{a} = (a(t,s); |t|+|s|
=1) \in[0,1)^4 $ is the (random) Green function defined as
%
\begin{eqnarray}
\label{green} g(t,s,\mathbf{a})&:=&\sum_{k=0}^\infty
a^{\star k}(t,s),
\end{eqnarray}
where $a^{\star k}(t,s)$ is the $k$-fold convolution of $a(t,s),
(t,s) \in\mathbb{Z}^2 $ defined recursively
by
\begin{eqnarray*}
a^{\star 0}(t,s) &=& \delta(t,s) := %
\cases{ 1, & \quad $(t,s)=
(0,0)$, \vspace*{3pt}
\cr
0, &\quad $(t,s) \ne(0,0)$,}
\\
a^{\star k }(t,s) &=& \sum_{(u,v) \in\mathbb{Z}^2}
a^{\star
(k-1)}(u,v) a(t-u,s-v), \qquad k \ge1.
\end{eqnarray*}
Note that \eqref{green}
can be rewritten as
%
\begin{equation}
\label{green0} g(t,s,\mathbf{a}) = \sum_{k=0}^\infty
A^k p_k(t,s),
\end{equation}
cf. \eqref{green00},
where $A$ is defined in \eqref{Acond} and $p_k(t,s) = \mathrm{P}(W_k
= (t,s)
|W_0 = (0,0))$ is the $k$-step probability of nearest-neighbor random
walk $\{ W_k, k =0,1,\ldots\} $ on $\mathbb{Z}^2 $
with one-step transition
probabilities
%
\begin{equation}
\label{transP} p(t,s) := \frac{a(t,s)}{A} \ge0, \qquad (t,s) \in
\mathbb{Z}^2.
\end{equation}
Generally, the $p_k(t,s)$'s
depend also on $\mathbf{a} = (a(t,s); |t|+|s| =1) \in[0,1)^4$ but this
dependence is suppressed below for brevity. Note that the series in
\eqref{green0} absolutely converges a.s., moreover,
%
\begin{equation}
\label{gsum} \sum_{(t,s) \in\mathbb{Z}^2} g(t,s,\mathbf{a}) = \sum
_{k=0}^\infty A^k \sum
_{(t,s) \in\mathbb{Z}^2} p_k(t,s) = \sum
_{k=0}^\infty A^k = \frac
{1}{1-A} <
\infty\qquad \mbox {a.s.}
\end{equation}
according to \eqref{Acond}. From \eqref{gsum}, it follows that
the Fourier transforms
$\hat p(x,y) := \sum_{|t|+|s|=1}\mathrm{e}^{-{\mathrm{i}}(tx+ sy)}
p(t,s)$ and
\begin{eqnarray*}
\hat g(x,y, \mathbf{a}) &:=& \sum_{(t,s) \in\mathbb{Z}^2} \mathrm
{e}^{-{\mathrm{i}
}(tx+ sy)} g(t,s,\mathbf{a}) = \sum_{k=0}^\infty
A^k \sum_{(t,s) \in\mathbb{Z}^2} \mathrm{e}^{-{\mathrm{i}
}(tx+ sy)}
p_k(t,s)
\\
&=& \sum_{k=0}^\infty A^k
\bigl(\hat p(x,y)\bigr)^k = \frac{1}{ 1 - A \hat p(x,y)}
\end{eqnarray*}
are well-defined and continuous on $\Pi^2:= [-\pi, \pi]^2$, a.s.
From Parseval's identity,
%
\begin{equation}
\label{parseval} \sum_{(t,s) \in\mathbb{Z}^2} \bigl|g(t,s,
\mathbf{a})\bigr|^2 = (2\pi)^{-2} \int_{\Pi^2}
\frac{ {\mathrm{d}} x \,\mathrm{d}y}{ |1 - A \hat p(x,y)|^2 }.
\end{equation}
Let
%
\begin{eqnarray}
q_1&:=&p(0,1) + p(0,-1) = 1 - p(1,0) - p(-1,0) =: 1-
q_2,\qquad q := \min(q_1, q_2),
\nonumber
\\[-8pt]
\label{qcoef}
\\[-8pt]
\nonumber
\mu_1&:=&p(1,0) - p(-1,0), \qquad \mu_2 :=
p(0,1) - p(0,-1),\qquad \mu := \sqrt{\mu_1^2 +
\mu_2^2}.
\end{eqnarray}
Note $q_i \in[0,1]$ and $q_1=0$ (resp., $q_2=0$) means that random
walk $\{ W_k \} $ is concentrated on the horizontal
(resp., vertical) axis of the lattice $\mathbb{Z}^2 $. Condition $\mu
= 0$
means that $\{ W_k \} $ has zero mean. Denote
%
\begin{equation}
\label{Psidef} \Psi(A,q,\mu) := \min \biggl( \frac{1}{q(1-A)},
\frac{1}{\mu
\sqrt
{q(1-A)}} \biggr) \biggl(1 + \log_+ \biggl( \frac{\mu^2}{q(1-A)}\biggr)
\biggr).
\end{equation}

The main result of this section is Theorem~\ref{exist}, below, which
provides sharp sufficient conditions for the convergence
of the series in \eqref{stationary0} involving the quantity $\Psi
(A,q,\mu)$ in \eqref{Psidef}. The proof of Theorem~\ref{exist} uses the
following Lemma~\ref{lemPsi}. The proof of this lemma
is given at the end of this subsection.

\begin{lem} \label{lemPsi} There exists a (non-random) constant
$C<\infty$ such that
%
\begin{equation}
\label{inPsi} \int_{\Pi^2} \frac{ \mathrm{d}x \,\mathrm{d}y}{|1 - A \hat
p(x,y)|^2} \le C \Psi (A,q,
\mu).
\end{equation}
\end{lem}

\begin{theorem} \label{exist} \textup{(i)} Assume there exists  $0< p \le2$  such that
%
\begin{equation}
\label{pmom} \mathrm{E}\bigl|\varepsilon (0,0)\bigr|^p < \infty \quad
\mbox{and} \quad \mathrm{E}\varepsilon (0,0) = 0 \qquad \mbox{for }1\le p \le2.
\end{equation}
Then there exists a stationary solution of random-coefficient equation
\eqref{NN} given by
\eqref{stationary0}, where the series converges conditionally a.s. and
in $L^p $ for any $\mathbf{a} = (a(t,s)\ge0, |t|+|s|=1) \in[0,1)^4 $ satisfying
\eqref{Acond}.

\textup{(ii)} In addition to \eqref{pmom}, assume that $q >0 $ a.s. and\vspace*{-2pt}
%
\begin{eqnarray}
\label{statcond} \cases{ \mathrm{E} \bigl[ \Psi(A,q,\mu)^{p-1}
(1-A)^{p-2} \bigr] < \infty, &\quad$\mbox{if } 1 < p \le2$,
\vspace*{3pt}
\cr
\mathrm{E} \bigl[ (1-A)^{2p-3} \bigr] < \infty, &\quad
$\mbox{if } 0 < p \le 1$.}
\end{eqnarray}
Then the series in \eqref{stationary0} converges unconditionally in
$L^p$, moreover,\vspace*{-2pt}
%
\begin{eqnarray}
\label{statmom} \mathrm{E} \bigl[ \bigl|X(t,s)\bigr|^p \bigr] &\le& C\cases{
\mathrm{E} \bigl[ \Psi(A,q,\mu)^{p-1} (1-A)^{p-2}
\bigr] < \infty, &\quad $1 < p \le 2$,\vspace*{3pt}
\cr
\mathrm{E} \bigl[
(1-A)^{2p-3} \bigr] < \infty, & \quad $0 < p \le 1$.}
\end{eqnarray}
\end{theorem}

\begin{pf}
Part (i) follows similarly as in \cite{ps2009},
proof of Proposition~1.
Let us prove part (ii). We shall use the following inequality; see
\cite{von1965}, also \cite{ps2009}, (2.7).
Let $0< p \le2$, and let $\xi_1, \xi_2, \ldots $ be random variables
with $\mathrm{E}|\xi_i|^p < \infty$. For $1\le p \le2$, assume in addition
that the $\xi_i$'s are independent and have zero mean
$\mathrm{E}\xi_i = 0$. Then
$ \mathrm{E} |\sum_i \xi_i  |^p  \le 2 \sum_i \mathrm
{E}|\xi_i|^p$.
The last inequality and the fact that~\eqref{stationary0} converges
conditionally in $L^p$
(see part (i))
imply that
%
\begin{equation}
\label{momp} \mathrm{E} \bigl[ \bigl|X(t,s)\bigr|^p | \mathbf{a} \bigr]
\le2\mathrm {E}\bigl|\varepsilon (0,0)\bigr|^p \sum
_{(u,v) \in\mathbb{Z}^2 } \bigl|g(u,v,\mathbf{a})\bigr|^p.
\end{equation}
Accordingly,
it suffices to prove\vspace*{-2pt} that
%
\begin{eqnarray}
\label{Lpbdd} \mathrm{E}\sum_{(t,s)\in\mathbb{Z}^2} \bigl|g(t,s,
\mathbf{a})\bigr|^p&<&\infty.
\end{eqnarray}
For $p=2$, \eqref{Lpbdd} is immediate from \eqref{parseval} and
\eqref
{inPsi}. Next,
using \eqref{parseval}, \eqref{inPsi} and H\" older's inequality,
for any $1< p < 2$ we obtain
%
\begin{eqnarray}
\sum_{(t,s) \in\mathbb{Z}^2} \bigl|g(t,s,\mathbf{a})\bigr|^p&=&
\sum_{(t,s) \in
\mathbb{Z}^2} \bigl|g(t,s,\mathbf{a})\bigr|^{2(p-1)}
\bigl|g(t,s,\mathbf{a})\bigr|^{2-p}
\nonumber
\\
\label{Lpbdd1} &\le& \biggl(\sum_{(t,s) \in\mathbb{Z}^2} \bigl|g(t,s,
\mathbf{a})\bigr|^2 \biggr)^{p-1} \biggl( \sum
_{(t,s) \in\mathbb{Z}^2} \bigl|g(t,s,\mathbf{a})\bigr| \biggr)^{2-p}
\\
&\le&C \Psi(A,q,\mu)^{p-1}(1-A)^{p-2}.
\nonumber
\end{eqnarray}
Next, consider the case $0< p \le1$. Using \eqref{green0}, the
inequality $|\sum_i x_i|^p \le\sum_i |x_i|^p $
and H\"older's inequality, we obtain
%
\begin{eqnarray}
\sum_{(t,s) \in\mathbb{Z}^2} \bigl|g(t,s,\mathbf{a})\bigr|^p&\le&
\sum_{k=0}^\infty A^{kp} \sum
_{|t|+|s| \le k} p^p_k(t,s)
\nonumber
\\[-2pt]
\label{Lpbdd2} &\le&\sum_{k=0}^\infty
A^{kp} \biggl\{ \sum_{|t|+|s| \le k}
p_k(t,s) \biggr\}^{p} \biggl\{ \sum
_{|t|+|s| \le k} 1 \biggr\}^{1-p}
\\[-2pt]
&\le&C\sum_{k=0}^\infty A^{kp}
k^{2(1-p)} \le \frac{C}{(1-
A^p)^{3-2p}} \le \frac{C}{(1- A)^{3-2p}},
\nonumber
\end{eqnarray}
where the last inequality follows from $1- x^p \ge p(1-x), x \in[0,1]$.
Note that $C$ in \eqref{Lpbdd1}--\eqref{Lpbdd2} are non-random. Hence,
\eqref{Lpbdd} follows from
\eqref{statcond} and the bounds in
\eqref{Lpbdd1}--\eqref{Lpbdd2}, proving the unconditional convergence of
\eqref{stationary0}.
Inequality \eqref{statmom}
is a consequence of \eqref{Lpbdd1}--\eqref{Lpbdd2} and~\eqref{momp}.
\end{pf}

\begin{pf*}{Proof of Lemma~\protect\ref{lemPsi}}
Write $I$ for the left-hand side
of \eqref{inPsi}. Since
\eqref{inPsi} holds trivially for $0\le A \le1/2$, we assume
$1/2 < A < 1$ in the sequel. We have
\begin{eqnarray*}
\label{Ap1} 1 - A \hat p(x,y) &=& (1-A) + A \sum
_{|t|+|s|=1} p(t,s) \bigl(1 - {\mathrm{e}}^{{\mathrm
{i}}(tx +
sy) }\bigr)
\\[-2pt]
&=& (1-A) + A \bigl[q_2\bigl(1- \cos(x)\bigr) + q_1
\bigl(1- \cos(y)\bigr) \bigr] - {\mathrm {i}} A \bigl(\mu _1 \sin(x) +
\mu_2 \sin(y)\bigr)
\end{eqnarray*}
and
\begin{eqnarray}
\bigl|1 - A \hat p(x,y)\bigr|^2 &=& \bigl((1-A) + A \bigl[q_2
\bigl(1- \cos(x)\bigr) + q_1 \bigl(1- \cos(y)\bigr) \bigr]
\bigr)^2
\nonumber
\\[-2pt]
&&{}+ A^2 \bigl(\mu_1 \sin(x) + \mu_2
\sin(y)\bigr)^2
\nonumber
\\[-9pt]
\label{Ap2}
\\[-9pt]
\nonumber
&\ge& (1/4) \bigl\{\bigl((1-A) + q \bigl[\bigl(1- \cos(x)\bigr) +
\bigl(1- \cos(y)\bigr) \bigr]\bigr)^2
\nonumber
\\[-2pt]
&&{}+ \mu^2 \bigl(\nu_1 \sin(x) + \nu_2
\sin(y)\bigr)^2 \bigr\},
\nonumber
\end{eqnarray}
where $\nu_i := \mu_i/\mu, i=1,2,  \nu_1^2 + \nu_2^2 = 1$.
Split $I = I_1 + I_2$, where $I_1 := \int_{[-\pi/4, \pi/4]^2}, I_2 :=
\int_{\Pi^2 \setminus[-\pi/4, \pi/4]^2}$.
Changing\vspace*{1pt} the coordinates $\sin(x) = u, \sin(y) = v,
\nu_1 u + \nu_2 v = s, -\nu_2 u + \nu_1 v = t, r^2 = t^2 + s^2, s = r
\sin(\phi)$ we get
\begin{eqnarray*}
I_1 &=&C \int_{[-1/\sqrt{2},1/\sqrt{2}]^2} \frac{1}{\sqrt
{(1-u^2)(1-v^2)}}
\\[-2pt]
&&{}\times\frac{\mathrm{d}u \,\mathrm{d}v}{  \{((1-A)+
q[(1-\sqrt{1-u^2})+
(1-\sqrt{1-v^2})])^2 + \mu^2 (\nu_1u + \nu_2 v)^2  \}}
\\[-2pt]
&\le&C\int_{u^2 + v^2 \le1} \frac{\mathrm{d}u\, \mathrm{d}v}{((1-A)
+ q[u^2 + v^2])^2
+ \mu^2 (\nu_1 u + \nu_2 v)^2}
\\[-2pt]
&=&C\int_{t^2 + s^2 \le1} \frac{\mathrm{d}s\, \mathrm{d}t} {
((1-A) + q [s^2 + t^2 ])^2 + \mu^2 s^2}
\\[-2pt]
&=&C \int_0^1\!\int_0^{\pi/2}
\frac{r \mathrm{d}r \,\mathrm{d}\phi
}{((1-A) + q
r^2)^2 +
\mu^2 r^2 \sin^2 (\phi)}.
\end{eqnarray*}
Using $\sin(\phi) \ge(1/2)\phi, \phi\in[0,\pi/2) $ with $W :=
\frac
{((1-A) + qx)^2 }{x} $, we obtain
\begin{eqnarray*}
I_1 &\le&C \int_0^1\! \int
_0^{1} \frac{\mathrm{d}x \,\mathrm{d}y}{((1-A)
+ q x)^2 + \mu
^2 x
y^2} \\
&\le & C\int
_0^1 \frac{\mathrm{d}x}{x} \int
_0^{1} \frac{\mathrm
{d}y}{ W + \mu^2
y^2}
\\
&\le&C\int_0^1 \frac{\mathrm{d}x}{x \sqrt{W} } \int
_0^{1/\sqrt
{W}} \frac
{\mathrm{d}
u}{1 + \mu^2 u^2}
\\
&\le&\frac{C}{\mu} \int_0^1
\frac{\mathrm{d}x}{x \sqrt{W} } \min\biggl(1, \frac{\mu
}{\sqrt{W}}\biggr) = C
\bigl(I'_1 + I''_1
\bigr),
\end{eqnarray*}
where
\begin{eqnarray*}
I'_1 &:=&\frac{1}{\mu} \int_0^1
\frac{\mathrm{d}x}{x \sqrt{W} } \mathbf{1} (\mu> \sqrt {W}), \qquad I''_1
:= \int_0^1 \frac{\mathrm{d}x}{x W } \mathbf{1} (
\mu< \sqrt{W}).
\end{eqnarray*}
Here,
%
\begin{eqnarray}
I_1'' &\le& \min \biggl( \int
_0^\infty\frac{\mathrm{d}x}{((1-A)+ q x)^2}, \frac
{1}{\mu}
\int_0^\infty\frac{\mathrm{d}x}{x^{1/2}((1-A) + q x)} \biggr)
\nonumber
\\[-8pt]
\label{I12}
\\[-8pt]
\nonumber
&\le& C \min \biggl(\frac{1}{q(1-A)}, \frac{1}{\mu\sqrt
{q(1-A)}} \biggr).
\end{eqnarray}
Since $I'_1 = 0$ for $\mu^2 \le q(1-A)$ we obtain
%
\begin{eqnarray}
I'_1 &\le & \frac{1}{\mu} \int
_0^1 \frac{\mathrm{d}x}{x^{1/2}
((1-A) + q x)
} \mathbf{1} \bigl(
\mu^2 > q(1-A)\bigr)
\nonumber
\\[-8pt]
\label{I11}
\\[-8pt]
\nonumber
&\le &  \frac{C}{\mu\sqrt{q(1-A)}} \mathbf{1} \bigl(
\mu^2 > q(1-A)\bigr).
\end{eqnarray}
Relations \eqref{I12} and \eqref{I11} yield
%
\begin{equation}
\label{I1} I_1 \le C \min \biggl( \frac{1}{q(1-A)},
\frac{1}{\mu\sqrt
{q(1-A)}} \biggr).
\end{equation}
Below we prove the bound
%
\begin{equation}
\label{I2} I_2 \le C %
\cases{ \displaystyle
(1-A+q)^{-2}, & \quad $\mu\le1-A+q$, \vspace*{3pt}
\cr
\displaystyle
\mu^{-1} (1-A+q)^{-1}\bigl(1+ \log\bigl(\mu/(1-A+q)\bigr)
\bigr), & \quad $\mu> 1-A+q$, }
\end{equation}
with $C$ independent of $A, q, \mu$, as elsewhere in this proof.
Since $1-A+q \ge\sqrt{q(1-A)}$, the desired inequality \eqref{inPsi},
viz., $I \le C \Psi(A, q,\mu)$,
follows from \eqref{I1} and \eqref{I2}.

Let us prove \eqref{I2}. For $\mu\le1 - A +q $ it follows trivially
from \eqref{Ap2}. Let $\mu> 1-A+q$ in the rest
of the proof.
From \eqref{Ap2}, we obtain that
%
\begin{eqnarray}
I_2 &\le&C \int_{\Pi^2 \setminus[-\pi/4,\pi/4]^2} \frac{\mathrm{d}x\,
\mathrm{d}y}{(1-A+
q)^2 + \mu^2 (\nu_1 \sin(x) + \nu_2 \sin(y))^2}
\nonumber
\\[-8pt]
\label{App}
\\[-8pt]
\nonumber
&\le&C\int_{[0,\pi/2]^2} \frac{\mathrm{d}x \,\mathrm{d}y}{(1-A +
q)^2 + \mu^2
(\tilde
\nu_1 \sin(x) + \tilde\nu_2 \sin(y))^2},
\end{eqnarray}
where $|\tilde\nu_i| = |\nu_i|, i=1,2 $ satisfy $\tilde\nu_1^2 +
\tilde\nu_2^2 = 1$.
Then
\begin{eqnarray*}
I_2&\le&C\int_{[0,1]^2} \frac{(1-u^2)^{-1/2} (1-v^2)^{-1/2}\, \mathrm
{d}u \,\mathrm{d}
v}{(1-A + q)^2 + \mu^2 (\tilde\nu_1 u + \tilde\nu_2 v)^2}
\\
&\le&\frac{C}{\mu^2}\int_{[0,1]^2} \frac{\mathrm{d}u \,\mathrm
{d}v}{(\epsilon^2 +
(\tilde\nu_1 u + \tilde\nu_2 v)^2)\sqrt{(1-u)(1-v)}}\qquad
\mbox{with } \epsilon:= \frac{1-A+q}{\mu} \ge0.
\end{eqnarray*}
We claim that
%
\begin{eqnarray}
\label{Jep} \int_{[0,1]^2} \frac{\mathrm{d}u \,\mathrm{d}v}{(\epsilon^2 +
(\tilde\nu_1 u +
\tilde
\nu_2 v)^2)\sqrt{(1-u)(1-v)}} &\le&
\frac{C}{\epsilon} \bigl(1 + \log_+ (1/\epsilon) \bigr)
\end{eqnarray}
with $C < \infty$ independent of $\epsilon>0$ and $\tilde\nu_i,
i=1,2, \tilde\nu_1^2 + \tilde\nu_2^2 = 1$.
Bound \eqref{Jep} proves \eqref{I2} and hence \eqref{inPsi} and the
lemma, too.
Therefore, it remains to prove \eqref{Jep}.

By symmetry, it suffices to prove \eqref{Jep} for $\tilde\nu_1 \ge
1/\sqrt{2}, 0 \ge\tilde\nu_2 \ge- 1/\sqrt{2}$, or
%
\begin{eqnarray}
J &:=& \int_{[0,1]^2} \frac{\mathrm{d}u \,\mathrm{d}v}{(\epsilon^2 + (u
- r v)^2)\sqrt
{(1-u)(1-v)}} \le
\frac{C}{\epsilon} \bigl(1 + \log_+ (1/\epsilon) \bigr) \nonumber
\\[-12pt]
\label{Jep1}
\\[-4pt]
\eqntext{\displaystyle
\mbox{uniformly in
} r \in[0,1].}
\end{eqnarray}
We have $J = \frac{1}{\sqrt{r}} \int_0^r \frac{\mathrm{d}v}{\sqrt{r-v}}
\int_0^1 \frac{\mathrm{d}u }{(\epsilon^2 + (u - v)^2)\sqrt{1-u}}
= \frac{1}{\sqrt{r}} \int_0^r \frac{\mathrm{d}v}{\sqrt{r-v}} \int_0^1
\frac{\mathrm{d}
z }{(\epsilon^2 + (1-z - v)^2)\sqrt{z}}
= \sum_{i,j=1}^2 J_{ij}$, where
\begin{eqnarray*}
J_{11} &:=&\frac{1}{\sqrt{r}} \int_0^r
\frac{\mathbf{1}(1-v > 2\epsilon)
\,\mathrm{d}
v}{\sqrt{r-v}} \int_0^1
\frac{ \mathbf{1}(|1-z-v|> \epsilon) \,\mathrm{d}z }{(\epsilon
^2 + (1-z -
v)^2)\sqrt{z}},
\\
J_{12} &:=&\frac{1}{\sqrt{r}} \int_0^r
\frac{\mathbf{1}(1-v > 2\epsilon)
\,\mathrm{d}
v}{\sqrt{r-v}} \int_0^1
\frac{ \mathbf{1}(|1-z-v|< \epsilon) \,\mathrm{d}z }{(\epsilon
^2 + (1-z -
v)^2)\sqrt{z}},
\\
J_{21} &:=&\frac{1}{\sqrt{r}} \int_0^r
\frac{\mathbf{1}(1-v < 2\epsilon)
\,\mathrm{d}
v}{\sqrt{r-v}} \int_0^1
\frac{ \mathbf{1}(|1-z-v|> \epsilon) \,\mathrm{d}z }{(\epsilon
^2 + (1-z -
v)^2)\sqrt{z}},
\\
J_{22} &:=&\frac{1}{\sqrt{r}} \int_0^r
\frac{\mathbf{1}(1-v < 2\epsilon)
\,\mathrm{d}
v}{\sqrt{r-v}} \int_0^1
\frac{ \mathbf{1}(|1-z-v|< \epsilon) \,\mathrm{d}z }{(\epsilon
^2 + (1-z -
v)^2)\sqrt{z}}.
\end{eqnarray*}
Bound \eqref{Jep1} will be proved for each $J_{ij}, i,j=1,2$.

\textit{Estimation of} $J_{22}$.
We have
%
\begin{eqnarray}
J_{22} &\le&\frac{1}{\sqrt{r}} \int_0^r
\frac{\mathbf{1}(1-v < 2\epsilon)
\,\mathrm{d}
v}{\sqrt
{r-v}} \frac{1}{\epsilon^2} \int_0^1
\frac{ \mathbf{1}(|1-z-v|< \epsilon) \,\mathrm{d}z }{\sqrt{z}}
\nonumber
\\
\label{Jep22} &\le&\frac{1}{\sqrt{r}} \int_0^r
\frac{\mathbf{1}(1-v < 2\epsilon)
\,\mathrm{d}
v}{\sqrt
{r-v}} \frac{1}{\epsilon^2} \int_0^{3\epsilon}
\frac{ \mathrm{d}z }{\sqrt{z}} \le \frac
{C}{\sqrt{r}
\epsilon^{3/2}} \int_0^r
\frac{\mathbf{1}(1-v < 2\epsilon)\, \mathrm
{d}v}{\sqrt
{r-v}}
\\
&\le&\frac{C}{\sqrt{r} \epsilon^{3/2}} \int_{r-2\epsilon}^r
\frac
{\mathrm{d}
v}{\sqrt{r-v}} \mathbf{1}(r > 1-2\epsilon) \le \frac{C}{\epsilon
}
\mathbf{1}(r>1-2\epsilon).
\nonumber
\end{eqnarray}

\textit{Estimation of} $J_{21}$.
We have
%
\begin{eqnarray}
J_{21} &\le&\int_0^1
\frac{\mathbf{1}(1-rv < 2\epsilon) \,\mathrm{d}v}{\sqrt{1-v}} \int_0^1
\frac{ \mathbf{1}(|1-z-rv|> \epsilon) \,\mathrm{d}z
}{(1-z-rv)^2 \sqrt{z}}
\nonumber
\\[-8pt]
\label{Jep21}
\\[-8pt]
\nonumber
&=&\int_0^{2\epsilon} \frac{\mathrm{d}x}{\sqrt{x}}
\int_0^1 \frac{\mathbf{1}(|x-z| > \epsilon) \,\mathrm{d}z }{ (x-z)^2
\sqrt{z}} \le
\frac{1}{\epsilon} \int_0^2
\frac{ \mathrm{d}x}{\sqrt{x}} \int_0^\infty
\frac{\mathbf{1}(|x-z| >1)\, \mathrm{d}z}{(x-z)^2 \sqrt
{z}} \le \frac
{C}{\epsilon}
\end{eqnarray}
since the last double integral converges.

\textit{Estimation of} $J_{12}$.
We have
%
\begin{eqnarray}
J_{12} &\le&\frac{1}{\sqrt{r}} \int_0^r
\frac{\mathbf{1}(1-v > 2\epsilon)
\,\mathrm{d}
v}{\sqrt
{r-v}} \frac{1}{\epsilon^2} \int_0^1
\frac{ \mathbf{1}(|(1-v)-z|< \epsilon) \,\mathrm{d}z }{\sqrt
{z}}
\nonumber
\\
&\le& \frac{C}{\sqrt{r}} \int_0^r
\frac{\mathbf{1}(1-v > 2\epsilon)
\,\mathrm{d}
v}{\sqrt
{r-v}} \frac{1}{\epsilon^2} \int_{1-v-\epsilon}^{1-v+\epsilon}
\frac{ \mathrm{d}z }{\sqrt{z}}
\nonumber
\\[-8pt]
\label{Jep12}
\\[-8pt]
\nonumber
&\le&\frac{C}{\sqrt{r} \epsilon^2} \int_0^r
\frac{\mathbf{1}(1-v >
2\epsilon)
\,\mathrm{d}v}{\sqrt{r-v}} (\sqrt{1-v+\epsilon} - \sqrt{1-v - \epsilon})
\nonumber
\\
&\le&\frac{C}{\sqrt{r} \epsilon^2} \int_0^r
\frac{\mathbf{1} (1-v>
2\epsilon)\,\mathrm{d}
v}{\sqrt{r-v}} \frac{\epsilon}{\sqrt{1-v}} \le \frac{C \log
(1/\epsilon)}{\epsilon}.\nonumber
\end{eqnarray}
Indeed, if $r \in[0,1/2]$ then $\frac{1}{\sqrt{r}} \int_0^r \frac
{\mathbf{1}
(1-v> 2\epsilon)\,\mathrm{d}v}{\sqrt{r-v} \sqrt{1-v}} \le\frac
{C}{\sqrt
{r}} \int_0^r
\frac{\mathrm{d}w}{\sqrt{w}} \le C$, and if
$r \in[1/2,1], \epsilon\le1/2 $ then with $z = w - (2\epsilon+ r -1)$
\begin{eqnarray*}
\frac{1}{\sqrt{r}} \int_0^r
\frac{\mathbf{1} (1-v> 2\epsilon)\,\mathrm
{d}v}{\sqrt{r-v}
\sqrt{1-v}} &\le&C\int_0^r
\frac{\mathbf{1} (w> 2\epsilon+r -1)\,\mathrm{d}w}{\sqrt
{w} \sqrt
{1-r+w}}
\\
&\le& C %
\cases{\displaystyle \int_0^1
\frac{ \mathrm{d}z}{\sqrt{z} \sqrt{z +
2\epsilon}}, & \quad $2\epsilon> 1-r$, \vspace*{3pt}
\cr
\displaystyle\int
_0^r \frac{\mathrm{d}w}{\sqrt{w} \sqrt{1-r+w}}, &\quad $1-r \ge 2
\epsilon$,} \le C\log(1/\epsilon).
\end{eqnarray*}

\textit{Estimation of} $J_{11}$.
We have
%
\begin{eqnarray}
J_{11}  &\le& \frac{1}{\sqrt{r}} \int_{1-r}^1
\frac{\mathbf{1}(w > 2\epsilon)
\,\mathrm{d}
w}{\sqrt
{w- (1-r)}} \int_0^1
\frac{ \mathbf{1}(|z-w|> \epsilon) \,\mathrm{d}z }{(z -
w)^2\sqrt{z}}
\nonumber
\\[-8pt]
\label{Jep11}
\\[-8pt]
\nonumber
&\le& \frac{C}{\epsilon\sqrt{r}} \int_{(1-r)/\epsilon
}^{1/\epsilon}
\frac{ L(w) \mathbf{1}(w > 2) \,\mathrm{d}w}{\sqrt{w- ({1-r})/{\epsilon}}},
\end{eqnarray}
where $L(w) := \int_0^{\infty} \frac{ \mathbf{1}(|z-w|> 1) \,\mathrm{d}z
}{(z -
w)^2\sqrt
{z}} \le Cw^{-1/2} $ for $w \ge1$. W.l.g.,
let $\epsilon\in(0,1/2]$. First, let $(1-r)/\epsilon< 1 $, then $r
\in(1/2,1]$ and
$w - \frac{1-r}{\epsilon} > w/2 $ for $w > 2 $. The above facts imply that
%
\begin{equation}
\label{Jep111} J_{11} \le \frac{C}{\epsilon} \int
_{1}^{1/\epsilon} \frac{ \mathrm
{d}w}{w} \le
\frac{C \log(1/\epsilon)}{\epsilon}\qquad \mbox{when } (1-r)/\epsilon< 1,
\end{equation}
with $C$ independent of $r, \epsilon$. Next, let $(1-r)/\epsilon\ge1
$ then from \eqref{Jep11}, $L(w) = O(w^{-1/2})$ and the change
of variables $w - \frac{1-r}{\epsilon} = \frac{1-r}{\epsilon} x$
we obtain
%
\begin{eqnarray}
\label{Jep112} J_{11} &\le&\frac{C}{\epsilon\sqrt{r}} \int
_{0}^{r/(1-r)}\frac{ \mathrm{d}
x}{\sqrt
{x} \sqrt{1+x}} \le
\frac{C}{\epsilon\sqrt{r}} %
\cases{ \bigl(r/(1-r)\bigr)^{1/2}, & \quad
$r \in[0,3/4]$, \vspace*{3pt}
\cr
\log\bigl(r/(1-r)\bigr), &\quad $r\in[3/4,1]$,}
\nonumber
\\[-8pt]
\\[-8pt]
\nonumber
&\le&\frac{C \log(1/\epsilon)}{\epsilon} \qquad \mbox{when } (1-r)/\epsilon\ge1.
\end{eqnarray}
Bounds \eqref{Jep111}, \eqref{Jep112} prove \eqref{Jep1} for $J_{11}$,
thereby completing
the proof of \eqref{Jep1}. Lemma~\ref{lemPsi} is proved.
\end{pf*}

\subsection{Aggregation of autoregressive RF}

\begin{defn}\label{Zdef}
Write $\varepsilon \in D(\alpha)$, $0<
\alpha\le
2$ if
\begin{longlist}[(ii)]
\item[(i)] $\alpha= 2$ and $\mathrm{E}\varepsilon = 0,  \sigma^2 := \mathrm{E}
\varepsilon
^2 <
\infty$, or
\item[(ii)] $0 < \alpha< 2 $ and there exist some constants
$c_{1}, c_{2} \ge0, c_{1} + c_{2} > 0$
such that
\[
\lim_{x \to\infty} x^\alpha\mathrm{P}(\varepsilon > x) =
c_{1} \quad \mbox{and} \quad \lim_{x \to-\infty}
|x|^{\alpha} \mathrm {P}(\varepsilon \le x) = c_{2};
\]
moreover, $\mathrm{E}\varepsilon = 0$ whenever $1 < \alpha< 2 $,
while for $\alpha
=1$ we assume that
the distribution of $\varepsilon $ is symmetric.
\end{longlist}
\end{defn}

\begin{rem} \label{fell}
Condition $\varepsilon \in D(\alpha$) implies
that the r.v. $\varepsilon $ belongs to the domain of normal
attraction of an
$\alpha$-stable law; in other words,
%
\begin{equation}
\label{Z} N^{-1/\alpha}\sum_{i=1}^N
\varepsilon _i \limd Z, \qquad N \to\infty,
\end{equation}
where $Z$ is an $\alpha$-stable r.v.; see \cite{fell1966}, pages 574--581.
The characteristic function of r.v. $Z$ in~\eqref{Z} is given by
%
\begin{equation}
\label{Zomega} \mathrm{E}\mathrm{e}^{{\mathrm{i}}\theta Z} = \mathrm{e}^{-|\theta|^\alpha\omega(\theta)},
\qquad \theta\in \mathbb{R},
\end{equation}
where $\omega(\theta)$ depends only on $\operatorname{sign}(\theta)$ and
$\alpha
, c_1, c_2, \sigma$
in Definition~\ref{Zdef}. See, for example, \cite{fell1966}, pages~574--581.
\end{rem}

Let $\{X_i(t,s)\},  i=1,2, \ldots$ be independent copies of \eqref
{stationary0} with i.i.d. innovations $\varepsilon (t,s) \in D(\alpha
), 0<
\alpha\le2$.
The aggregated field $ \{ \mathfrak{X}(t,s); (t,s) \in\mathbb{Z}^2 \}
$ is
defined as the limit in distribution:
%
\begin{eqnarray}
\label{aggreNN} N^{-1/\alpha} \sum_{i=1}^N
X_i(t,s)&\limfdd&\mathfrak{X}(t,s), \qquad (t,s) \in
\mathbb{Z}^2, N \to\infty.
\end{eqnarray}
Introduce an independently scattered $\alpha$-stable random measure $M$
on $\mathbb{Z}^2 \times[0,1)^4 $
with characteristic functional
%
\begin{equation}
\label{Mdef} \mathrm{E}\exp \biggl\{{\mathrm{i}}\sum
_{(t,s) \in\mathbb{Z}^2} \theta_{t,s} M_{t,s}(B_{t,s})
\biggr\} = \exp \biggl\{ - \sum_{(t,s)\in\mathbb{Z}^2} |
\theta_{t,s}|^\alpha \omega (\theta_{t,s})
\Phi(B_{t,s}) \biggr\},
\end{equation}
where $\Phi(B) := \mathrm{P}( \mathbf{a} = (a(t,s), |t|+|s| =1) \in B)$
is the
mixing distribution,
$\theta_{t,s} \in\mathbb{R}$, $B, B_{t,s} \subset[0,1)^4 $ are arbitrary
Borel sets, and $\omega$ is the same as in
\eqref{Zomega}. According to the terminology in \cite{sam1994},
Definition~3.3.1, $M$ is called an $\alpha$-stable measure
with control measure $\operatorname{Re} (\omega(1)) \Phi(\mathrm{d}\mathbf{a}) $
proportional to the mixing distribution $\Phi$, and a constant
skewness intensity $\operatorname{Im}(\omega(1))/\operatorname{Re}(\omega(1))\tan(\pi \alpha/2)$.

\begin{proposition} \label{existmix} Let $\varepsilon (0,0) \in
D(\alpha),
0< \alpha\le2 $. Assume that the mixing distribution satisfies
the following condition: there exists $\epsilon>0$ such that
%
\begin{eqnarray}
\label{aggrecond} \cases{ \mathrm{E}\bigl[\Psi(A,q,\mu)\bigr] < \infty, & \quad $
\mbox{if } \alpha= 2$, \vspace*{3pt}
\cr
\mathrm{E}\bigl[\Psi^{p-1}(A,q,
\mu) (1-A)^{p-2}\bigr] < \infty, &\quad $\mbox{if } 1 < \alpha< 2, p =
\alpha\pm\epsilon$, \vspace*{3pt}
\cr
\mathrm{E}\bigl[(1-A)^{2\alpha-3 - \epsilon}\bigr] <
\infty, & \quad$\mbox {if } 0< \alpha \le1$, }
\end{eqnarray}
where $\Psi(A,q,\mu)$ is defined in \eqref{Psidef}.
Then the limit aggregated RF in \eqref{aggreNN} exists and has the
stochastic integral representation
%
\begin{eqnarray}
\label{aggremix0} \mathfrak{X}(t,s) &=&\sum_{(u,v) \in\mathbb{Z}^2}
\int_{[0,1)^4} g(t-u,s-v, \mathbf{a}) M_{u,v}(\mathrm{d}
\mathbf{a}),\qquad (t,s) \in\mathbb{Z}^2.
\end{eqnarray}
\end{proposition}

\begin{rem} \label{alpha}
Note for the 3N and 4N models, we have
$\mu= 1, q = 1/3, \Psi(A, 1/3, 1) \le\frac{C}{\sqrt{1-A}} (1+|\log(1-A)|)
$ and $\mu=0, q = 1/4, \Psi(A, 1/4, 0) \le C/(1-A)$, respectively. As
a consequence,
for the aggregated 3N and 4N models and a regularly varying (mixing)
density of~$A$ in \eqref{mixdensity}, condition \eqref{aggrecond} for $1 <
\alpha
\le2$ reduces to
$\beta> - (\alpha-1)/2 $ and $\beta>0$, respectively.
\end{rem}

\begin{pf*}{Proof of Proposition \protect\ref{existmix}}
Let $T \subset
\mathbb{Z}
^2 $ be a finite set, $\theta_{t,s} \in\mathbb{R},  (t,s) \in T$,
and $S_N =
N^{-1/\alpha} \sum_{i=1}^N U_i $ be a sum of i.i.d. r.v.'s with common
distribution
\begin{eqnarray*}
U &:=& \sum_{(t,s) \in T} \theta_{t,s} X(t,s) =
\sum_{(u,v) \in
\mathbb{Z}
^2} G(u,v,\mathbf{a}) \varepsilon (u,v),
\\
G(u,v,\mathbf{a}) &:= &\sum_{(t,s) \in T}
\theta_{t,s} g(t-u,s-v,\mathbf{a}).
\end{eqnarray*}
It suffices to prove that
$S_N \limd S  (N \to\infty)$, where
$S := \sum_{(t,s) \in T} \theta_{t,s} \mathfrak{X}(t,s)$ is a
$\alpha$-stable r.v. with characteristic function
\[
\mathrm{E}\mathrm{e}^{{\mathrm{i}} w S} = \exp \biggl\{ - |w|^\alpha \sum
_{(u,v) \in
\mathbb{Z}^2} \mathrm{E} \bigl[ \bigl| G(u,v,\mathbf{a})
\bigr|^\alpha \omega \bigl(w G(u,v,\mathbf{a}) \bigr) \bigr] \biggr\}.
\]
For this, it suffices to prove that r.v. $U$ belongs to the domain
of attraction of r.v. $S$ (in the sense of (\ref{Z})) or $U \in
D(\alpha
)$, see Remark~\ref{fell};
in other words, that
%
\begin{equation}
\label{U2} \mathrm{E}U^2 = \mathrm{E}S^2 < \infty
\qquad \mbox{for } \alpha= 2,
\end{equation}
and, for $0< \alpha< 2$,
%
\begin{eqnarray}
\lim_{x \to\infty} x^{\alpha} \mathrm{P}(U > x)&=& \sum
_{(u,v) \in
\mathbb{Z}^2} \mathrm{E} \bigl[ \bigl| G(u,v,\mathbf{a})
\bigr|^\alpha \bigl\{c_1 \mathbf{1}\bigl(G(u,v,\mathbf{a}) > 0
\bigr) + c_2 \mathbf{1}\bigl(G(u,v,\mathbf{a}) <0\bigr) \bigr\}
\bigr],
\nonumber
\\[-8pt]
\label{Ualpha}
\\[-5pt]
\nonumber
\lim_{x \to-\infty} |x|^{\alpha} \mathrm{P}(U \le
x)&=& \sum_{(u,v)
\in\mathbb{Z}
^2} \mathrm{E} \bigl[ \bigl| G(u,v,
\mathbf{a}) \bigr|^\alpha \bigl\{c_1 \mathbf{1}\bigl(G(u,v,
\mathbf{a}) < 0\bigr) + c_2 \mathbf{1}\bigl(G(u,v,\mathbf{a}) >0\bigr)
\bigr\} \bigr],
\end{eqnarray}
where $c_i, i=1,2 $ are the asymptotic constants in Definition~\ref
{Zdef} satisfied by $\varepsilon (0,0) \sim D(\alpha)$.
Here, \eqref{U2} follows from definitions of $U$ and $S$ and Theorem~\ref{exist} with $p=2$.
To prove \eqref{Ualpha},
we use \cite{hult2008}, Theorem~3.1. Accordingly, it suffices to
check that
there exists $\epsilon>0$ such that for $0< \alpha< 2,  \alpha\ne1$,
%
\begin{eqnarray}
\label{UG} &&\sum_{(u,v) \in\mathbb{Z}^2} \mathrm{E} \bigl| G(u,v,
\mathbf{a}) \bigr|^{\alpha+ \epsilon} < \infty\quad \mbox{and}\quad \sum
_{(u,v) \in\mathbb{Z}^2} \mathrm{E} \bigl| G(u,v,\mathbf{a}) \bigr|^{\alpha-
\epsilon
} <
\infty,
\end{eqnarray}
and
\[
\mathrm{E} \biggl(\sum_{(u,v) \in
\mathbb{Z}^2} \bigl| G(u,v,\mathbf{a})
\bigr|^{\alpha- \epsilon} \biggr)^{({\alpha+\epsilon})/({\alpha- \epsilon})} < \infty\qquad \mbox{for } \alpha= 1.
\]
Since $T\subset\mathbb{Z}^2$ is a finite set, it suffices to show
\eqref{UG}
with $G(u,v,\mathbf{a})$ replaced by $g(u,v,\mathbf{a})$.
Let $1< \alpha< 2$ and $p = \alpha\pm\epsilon\in(1, 2)$ in \eqref
{aggrecond}. Then
$\sum_{(u,v) \in\mathbb{Z}^2} \mathrm{E} |g(u,v,\mathbf{a})
|^{p} \le C\mathrm{E}
[\Psi
(A,q,\mu)^{p-1} (1-A)^{p-2}] < \infty$
follows from \eqref{Lpbdd1} and
\eqref{aggrecond}.
In the case $0< \alpha< 1$, relations \eqref{UG} immediately
follow from \eqref{Lpbdd2} and \eqref{aggrecond} with $p = \alpha\pm
\epsilon\in(0,1)$. For
$\alpha=1$, \eqref{UG} follows from \eqref{Lpbdd2} in a similar way.
\end{pf*}

\section*{Acknowledgements}

The authors are grateful to several anonymous referees for numerous
suggestions and comments
which helped to improve this and the previous versions of this paper.
This research
was supported by a grant (No. MIP-063/2013) from the Research Council
of Lithuania.






\printhistory
\end{document}